\documentclass[11pt,leqno,amscd,amssymb,verbatim, url]{amsart}
\newtheorem{thm}{Theorem}[section]
\usepackage{amsfonts,latexsym}
\newtheorem{lem}[thm]{Lemma}
\newtheorem{cor}[thm]{Corollary}

\newtheorem{prop}[thm]{Proposition}
\theoremstyle{definition}

\newtheorem{rem}[thm]{Remark}
\newtheorem{rems}[thm]{Remarks}

\numberwithin{equation}{thm}

\newcommand{\sH}{{\mathcal H}}

\newcommand{\wSp}{{S}_p}

\newcommand{\sC}{{\mathcal {C}}}

\newcommand{\Ext}{{\text{\rm Ext}}}

\newcommand{\x}{\text{\bf \it x}}

\newcommand{\sS}{{\mathcal S}}

\newcommand{\Hom}{\text{\rm Hom}}

\newcommand{\End}{\operatorname{End}}

\newcommand{\ch}{\operatorname{ch}}

\newcommand{\Jan}{\Gamma_{\text{\rm Jan}}}

\newcommand{\uomega}{{\underline\omega}}
\newcommand{\ulambda}{{\underline\lambda}}
\newcommand{\utau}{{\underline\tau}}
\newcommand{\umu}{{\underline\mu}}
\newcommand{\unu}{{\underline\nu}}
\newcommand{\usigma}{{\underline\sigma}}
\newcommand{\ugamma}{{\underline\gamma}}

\newcommand{\blist}{\begin{list}{\rom{(\roman{enumi})}}{\setlength{\leftmarg
in}{0em}

\setlength{\itemindent}{7ex}
\setlength{\labelsep}{2ex}\setlength{\listparindent}{\parindent}
\usecounter{enumi}}}
\newcommand{\elist}{\end{list}}

\dedicatory{We dedicate this paper to the memory of Walter Feit}
\begin{document}

\title[Extensions, Levi subgroups and character formulas]{\large {\bf Extensions,
 Levi subgroups and character formulas}}
\author{Brian J. Parshall}
\address{Department of Mathematics \\
University of Virginia\\
Charlottesville, VA 22903} \email{bjp8w@virginia.edu {\text{\rm
(Parshall)}}}
\author{Leonard L. Scott}
 \address{Department of Mathematics \\
University of Virginia\\
Charlottesville, VA 22903} \email{lls2l@virginia.edu {\text{\rm
(Scott)}}}
 \begin{abstract} This paper consists of three
interconnected parts. In Parts I,III we study the relationship
between the cohomology of a reductive group $G$ and that of a Levi
subgroup $H$.  For example, we provide a necessary condition,
arising from Kazhdan-Lusztig theory,  for a natural map
$\Ext^\bullet_G(L,L')\to\Ext^\bullet_H(L_H,L'_H)$ to be surjective,
given irreducible $G$-modules $L,L'$ and corresponding irreducible
$H$-modules $L_H,L'_H$. In cohomological degree $n=1$, the map is
always an isomorphism, under our hypothesis.   These results were
inspired by recent work of Hemmer \cite{H} obtained for $G=GL_n$,
and both extend and improve upon the latter when our condition is
met.
 Part II
obtains results on Lusztig character formulas for reductive groups,
obtaining new necessary and sufficient conditions for such formulas
to hold.  In the special case of $G=GL_n$, these conditions can be
recast in a striking way completely in terms of explicit
representation theoretic properties of the symmetric group (and the
results  improve upon the sufficient cohomological conditions
established recently in \cite{PS}).
\end{abstract}
\thanks{Research supported in part by the National Science
Foundation}


\maketitle {\SMALL \tableofcontents}
\section*{Introduction}
Let $H$ be a Levi subgroup of a reductive  algebraic group $G$
defined over an algebraically closed field $k$ of positive
characteristic $p$. This paper has its origins in the authors'
attempts to understand the relationship between the cohomology of
$H$ and that of $G$. Usually, we concentrate on comparing
$\Ext^\bullet_G(M,N)$ to $\Ext^\bullet_H(M_H,N_H)$, where $M,N$
are rational $G$-modules and $M_H,N_H$ are rational $H$-modules
obtained from $M,N$ by a natural ``truncation" process. Especially
important examples include the cases when $M,N$ are taken to be
standard (or Weyl) modules, costandard modules, or irreducible
modules.

Contributions to the problem have been made by E. Cline \cite{C}
and S. Donkin \cite{D1}.  Their work was the starting point for
our discovery in \cite{CPS6} (with E. Cline) of an explicit
category equivalence between some highest weight categories
associated to the categories of rational modules for $G$ and $H$.
Because the methods were generally formal, they also applied to
module categories for quantum enveloping algebras and, in type
$A$, for $q$-Schur algebras at a root of unity. Furthermore, these
results provided natural $\Ext$-transfer results in going from
certain $H$-modules to $G$-modules.

 Part I of the present paper reviews the category equivalence described above,
 and gives a number of examples not explicitly considered in \cite{CPS6}.
 For example, when $H$ decomposes into a direct product of reductive subgroups,
 the groups $\Ext_H^\bullet(M_H,N_H)$ can be further analyized by means of an
 evident application of the K\"unneth theorem.    In type $A$,
these results can be stated in an elegant way, using the
elementary combinatorics of partitions. In this case, our
formulation was inspired by the corresponding degree $0$
decomposition given by Lyle-Mathas in a preprint \cite{LM}. (More
recently, A. Mathas has made a revised verison of \cite{LM}
available to us, in which they also obtain, by different methods,
the K\"unneth decomposition for $q$-Schur algebras.) It should be
remarked, however, that the methods of \cite{CPS6} apply for
reductive groups or quantum groups in {\it all types}. As an
illustration, we very briefly consider the case of type $C$, which
behaves like the type $A$ case.

The results of Part II improve upon \cite[\S7]{PS} which
formulated a vanishing condition for symmetric group cohomology
which, if true, would be sufficient to prove the Lusztig character
conjecture for $G=SL_n(k)$. In
\S\ref{equivalentconditionssection}, we allow $G$ to be of any
type with $p\geq h$, the Coxeter number of $G$. Roughly, we
consider finite ideals $\Gamma$ of dominant weights having the
following property ($\dagger$): if $\varpi\in\Gamma$ is regular,
then the character of the irreducible $G$-module $L(\varpi)$ of
highest weight $\varpi$ is given by the Lusztig character formula;
see (\ref{lusztigcharacterformula}). Lusztig's conjecture \cite{L}
asserts that the Jantzen region $\Jan$ (defined in
(\ref{jantzenregion})) is such an ideal; however, independently of
that conjecture, finite ideals $\Gamma$ satisfying property
($\dagger$) are easily shown to exist. In that  case, various
homological methods can be brought to bear on the category
${\mathcal C}_G[\Gamma]$ of finite dimensional rational
$G$-modules which have composition factors $L(\gamma)$,
$\gamma\in\Gamma$.
 For
example, \cite{CPS2} establishes that  $\Gamma$ satisfies
($\dagger$) if and only if $\Ext^1_G(L(\lambda),L(\mu))\not=0$
whenever $\lambda,\mu\in\Gamma$ are regular weights which are
mirror images of each other in adjacent alcoves. This result is
recalled in Thm. \ref{lusztigequiv}, which also presents a useful
(and new) method allowing the possible enlargement of (certain
ideals) $\Gamma$ by the addition of an alcove adjacent to
$\Gamma$. Thm. \ref{firstreduction} recasts the $\Ext^1$-criterion
above
 in terms of a condition on mappings between
$\Ext^1$-groups. The modules involved are no longer irreducible;
indeed, all are standard modules or modules with a standard
filtration.

 In Section \ref{typeA}, Thm. \ref{typeareduction} translates, in the
case $G=SL_n(k)$, the results of
\S\ref{equivalentconditionssection} into equivalent results stated
{\it entirely} in symmetric group terms. In particular, we obtain
a necessary and sufficient condition for the Lusztig character
conjecture to hold for $SL_n(k)$, stated entirely in terms of the
representation theory of the symmetric group. This result is in
contrast to \cite[\S7]{PS}, which only obtained a sufficient
condition. Section \ref{typeA} concludes with a discussion of the
relation of Thm. \ref{typeareduction} to the sufficient condition
established in \cite[\S7]{PS}. See Remark \ref{erdmannremark}(b)
for a brief comparison between features of our Lusztig
 conjecture reduction to those of Erdmann's general computation \cite{E2} of
 Schur algebra decomposition numbers in terms of  those of the symmetric group.

The results in Part I for $q$-Schur algebras apply to
$\Ext$-groups between standard modules, or between costandard
modules, or between an irreducible and a standard or costandard
module. Part III takes up  the question of comparing $\Ext$-groups
between {\it irreducible} $G$-modules with those $\Ext$-groups
between similarly parametrized irreducible $H$-modules. One result
in this direction, for Schur algebras, has already been obtained
by Hemmer \cite{H}, who proves the injectivity of the natural map
$\phi_H^G({\varpi,\varpi'}):
\Ext^1_G(L(\varpi),L(\varpi'))\to\Ext^1_H(L_H(\varpi),L_H(\varpi'))$
when the dominant weights $\varpi,\varpi'$ lie in the same coset
relative to the weight lattice of $H$. His result is stated in the
type $A$ formalism of row removal in Young diagrams and Schur
algebras. However, in \S\ref{hemmersection}, Thm. \ref{Hemmer}
proves that Hemmer's result holds for all types.
 (A similar result would hold, with the same
proof, for quantum enveloping algebras.) Hemmer has asked whether
the map $\phi_H^G({\varpi,\varpi'})$ is an isomorphism. In
\S\ref{example}, we provide several (related) examples showing
that the answer to Hemmer's question is sometimes negative.  The
examples involve weights at which the Lusztig character formula
fail, and thus suggest a connection between Hemmer's question and
the validity of the Lusztig character formula.

Now assume, for general $G$, that the Lusztig character formula
holds for all irreducible modules $L(\varpi)$ with regular highest
weights $\varpi$ in a given finite ideal $\Gamma$ in the poset of
dominant weights. Then, as perhaps suggested by the example in
\S\ref{example}, Thm. \ref{kltheorem} shows that
$\phi^G_H(\varpi,\varpi')$ {\it is} an isomorphism if
$\varpi,\varpi'\in\Gamma$ are regular and lie in the same coset
relative to the root lattice of $H$. Thus, Hemmer's question does
  have a positive answer in this case. Moreover, we prove
that the restriction map
$\Ext^n_G(L(\varpi),L(\varpi'))\to\Ext^n_H(L(\varpi),L_H(\varpi'))$
is {\it surjective in all degrees} $n$. A key ingredient in these
results is the study in \cite{CPS3} of the structure of the
algebra $A^!=\Ext^\bullet_A(L_0,L_0)$, when $A$ is a
quasi-hereditary algebra such that $A$--mod has a Kazhdan-Lusztig
theory; here $L_0$ is the direct sum of the distinct irreducible
$A$-modules.

Actually, Part III of this paper developed before Part II,
following quite naturally upon the investigation in Part I. The
relationship with the Lusztig character formula which emerged
inspired the revisiting of the latter topic in Part II, though
only the results up through Thm. \ref{lusztigequiv} are really
required in Part III.  One interesting by-product of this
investigation is a final
 observation, given in Cor. \ref{lastcor}, showing that  a form
 of the Lusztig conjecture holds for any Levi subgroup $H$, once it
 holds for $G$.

\medskip\medskip
\begin{center}{\large \bf PART I:  A category equivalence and applications}\end{center}
 \addcontentsline{toc}{section}{\bf Part I: A category equivalence and applications }
\medskip

 We first
review an essential category equivalence involving Levi subgroups
of reductive groups.  The method has applicability to other
situations, e.~g., quantum groups, $q$-Schur algebras, as we
illustrate. In addition, the equivalence will play an important
role later in the paper.

\section{The equivalence}\label{sectionone} Fix a reductive group $G$
over a fixed algebraically closed field $k$ of positive
characteristic $p$.  We will assume that the derived group of $G$
is simply connected. Let $B$ be a Borel subgroup of $G$ and let
$T$ be a maximal torus of $G$ contained in $B$. Let $X=X(T)$ be
the character group of $T$, and set $X^\vee=\Hom({\mathbb
G}_m,T)$, the cocharacter group of $T$. There is a natural perfect
pairing $\langle\,,\,\rangle:X\times X^\vee\to{\mathbb
Z}\cong\End({\mathbb G}_\alpha)$. We let $\Phi\subset X$ be the
root system of $T$ in $G$, and $\Phi^+$ (resp.,  $\Pi$) be the set
of positive (resp., simple) roots determined by $B$; thus, if
$\alpha\in\Phi^+$, the corresponding one-parameter root subgroup
$U_\alpha$ is a subgroup of $B$. For $\alpha\in\Phi$, let
$\alpha^\vee\in X^\vee$ be the associated coroot. Now let
$P\supseteq B$ be a parabolic subgroup, having Levi factor
$H\supseteq T$.  We also refer to $H$ as a {\it Levi subgroup}.
The root system of $T$ in $H$ is denoted $\Phi_H$, while
$\Phi^+_H:=\Phi_H\cap\Phi^+$ is a set of positive roots and
$\Pi_H=\Pi\cap\Phi_H$ is the corresponding set of simple roots.

We let $X^+\subset X$ (resp., $X^{+(H)}\subset X$) be the set of
dominant (resp., $H$-dominant) weights on $T$. Thus, $\lambda\in
X^+$ (resp., $\lambda\in X^{(H)+}$) provided that $\langle
\lambda,\alpha^\vee\rangle\in{\mathbb Z}^+$ for all $\alpha\in\Pi$
(resp., $\alpha\in\Pi_H$).  Define poset structures $\leq$ and
$\leq_H$ on $X$ by putting $\lambda\leq\mu$ (resp.,
$\lambda\leq_H\mu$) provided that $\mu-\lambda\in{\mathbb
Z}^+\Phi^+$ (resp., $\mu-\lambda\in{\mathbb Z}^+\Phi^+_H$). By
\cite[Lemma 6]{CPS6}, $X^+$ is an ideal in the poset
$(X^{(H)+},\leq_H)$, i.~e., if $\mu\leq_H\lambda$ and $\mu\in
X^{(H)+}$, $\lambda\in X^+$, then $\mu\in X^+$.

Let ${\mathcal C}_G$ (resp., ${\mathcal C}_H$) be the category of
finite dimensional rational $G$-modules (resp., $H$-modules). The
irreducible objects in $\mathcal C$ (resp., ${\mathcal C}_H$) are
indexed by the set $X^+$ (resp., $X^{(H)+}$).
   For $\lambda\in X^+$, let $L(\lambda)$,
$\Delta(\lambda)$, and $\nabla(\lambda)$ be the (rational)
irreducible, standard, and costandard $G$-modules, respectively,
of highest weight $\lambda$. Similarly, for $\sigma\in X^{(H)+}$,
$L_H(\sigma),\Delta_H(\sigma),\nabla_H(\sigma)$ be the analogous
objects in the category of rational $H$-modules.

Let $\Gamma^+$ be a finite ideal in the poset $(X^+,\leq)$, and
let ${\mathcal C}_G[\Gamma^+]$ be the full subcategory of
$\mathcal C$ consisting of all objects in $\mathcal C$ which have
composition factors $L(\gamma)$, $\gamma\in\Gamma^+$. Thus,
${\mathcal C}_G[\Gamma^+]$ is a highest weight category with poset
$(\Gamma^+,\leq)$.\footnote{If $\Gamma\subset X$, it is convenient
to write $\sC_G[\Gamma]$ for the full subcategory of $\sC_G$
having objects which have composition factors $L(\gamma)$,
$\gamma\in\Gamma$. Thus, if $\Gamma^+:=\Gamma\cap X^+$, then
$\sC_G[\Gamma]=\sC_G[\Gamma^+]$.  Also, if $\Gamma$ is an ideal in
$(X,\leq)$, $\Gamma^+$ is an ideal in $(X^+,\leq)$. \par The set
$X^+$ in this paragraph may be replaced with the set of highest
weights of irreducible modules in any block or union of blocks for
the category of rational $G$-modules. The same discussion applies.
We will return to this topic later in
\S\ref{extoneisomorhismandhigherext}. If will be useful to observe
that if $\Gamma$ is a union of blocks in an ideal, then ${\mathcal
C}_G[\Gamma]$ is a highest weight category with poset
$(\Gamma,\leq)$.} If $\Xi^+\subseteq\Gamma^+$ is a coideal (i.~e.,
$\Gamma^+\backslash\Xi^+$ is an ideal), we put ${\mathcal
C}_G(\Xi^+)$ equal to the quotient category of ${\mathcal
C}_G[\Gamma^+]$ by the Serre subcategory ${\mathcal
C}_G[\Gamma^+\backslash\Xi^+]$. Then ${\mathcal C}_G(\Xi^+)$ is a
highest weight category with poset $(\Xi^+,\leq)$.\footnote{The
definition of $\sC_G(\Xi^+)$ does not depend on the choice of the
finite ideal $\Gamma^+$: Suppose that $\Xi^+$ is any finite subset
of $X^+$ which is interval closed, in the sense that, given,
$\xi_1<\xi_2$ in $\Xi^+$, any $\xi\in X^+$ satisfying
$\xi_1\leq\xi\leq \xi_2$ automatically belongs to $\Xi^+$. Then
$\Xi^+$ is a coideal in the finite ideal $\Gamma'$ consisting of
all $\gamma\in X^+$ satisfying $\gamma\leq\omega$ for some
$\xi\in\Xi$. Thus, $\sC_G(\Xi^+)$ can be taken to be the highest
weight category $\sC_G[\Gamma']/\sC_G[\Gamma'\backslash\Xi^+]$. If
$\Xi^+$ is a coideal in another finite ideal $\Gamma^+$, it is
readily verified that $\sC_G(\Xi^+)$ is equivalent to
$\sC_F[\Gamma^+]/\sC_G[\Gamma^+\backslash\Xi^+]$.} Let
$j^*:{\mathcal C}_G[\Gamma^+]\to{\mathcal C}_G(\Xi^+)$ be the
natural quotient functor. For $\xi\in\Xi^+$, the functor $j^*$
takes the irreducible (resp., standard, costandard) object
$L(\xi)$ (resp., $\Delta(\xi)$,$\nabla(\xi)$) in $\sC_G[\Gamma^+]$
indexed by $\xi$ to the irreducible (resp., standard, costandard)
object $L_{\Xi^+}(\xi)$ (resp., $\Delta_{\Xi^+}(\xi)$,
$\nabla_{\Xi^+}(\xi)$) in $\sC_G(\Xi^+)$ indexed by $\xi$. If
$\gamma\in \Gamma^+\backslash\Xi^+$, $j^*$ annihilates
$L(\gamma)$, $\Delta(\gamma)$, $\nabla(\gamma)$. These facts are
all well-known general properties of highest weight
categories---see, e.~g., \cite[\S2]{CPS6} and the references given
there.

    Fix $\omega\in X$ and
let $\Omega=\omega+{\mathbb Z}\Phi_H$ be the corresponding coset
of ${\mathbb Z}\Phi_H$ in $X$.  Set $\Omega^+=\Omega\cap X^+$ and
$\Omega^{(H)+}=\Omega\cap X^{(H)+}$. Let $F$ be a finite, nonempty
subset of $\Omega^+$.\footnote{If $G$ has simple derived group
$G'$, then $\Omega^+$ is itself finite \cite[Prop. 9]{CPS6}, so we
can usually take $F=\Omega^+$. However, even in that case, it
might be desirable to consider a proper subset $F$.} Let
$\Gamma^+_F$ denote the (finite) ideal in $(X^+,\leq)$ generated
by $F$, and put $\Omega_F^+=\Gamma^+_F\cap \Omega$.
   It is shown in \cite[p. 222]{CPS6} that
$\Omega_F^+$ is a {\it coideal} in $(\Gamma_F^+,\leq)$, but an
{\it ideal} in $(X^{(H)+},\leq_H)$.  We consider the truncation
functor
\begin{equation}\label{truncationfunctor}
\pi_{\Omega }:{\mathcal C}_G[\Gamma^+_F]\to{\mathcal
C}_H[\Omega]={\mathcal C}_H[\Omega^{+(H)}],\quad M\mapsto
\pi_{\Omega}M=\bigoplus_{\tau\in\Omega}M_\tau.
\end{equation}
   Here $M_\tau$,
$\tau\in\Omega$, denotes the $\tau$-weight space of the rational
$G$-module $M$. By \cite{D1} and \cite[Prop. 7]{CPS6},
$\pi_\Omega$ maps $L(\omega),\Delta(\omega),\nabla(\omega)$,
$\omega\in\Omega^+_F$ to the corresponding objects $L_H(\omega),
\Delta_H(\omega),\nabla_H(\omega)$ in ${\mathcal
C}_H[\Omega^+_F]$. (If $\gamma\in\Gamma^+_F\backslash\Omega^+_F$,
$\pi_\Omega$ maps these objects to $0$.) We have:

\begin{thm}\label{maintheorem} \cite[Thm. 8]{CPS6} The
functor $\pi_{\Omega}$ factors through the quotient morphism
$j^*:{\mathcal C}_G[\Gamma_F^+]\to{\mathcal C}_G(\Omega^+_F)$ to
produce an equivalence
\begin{equation}\label{equiv1}
{\sC}_G(\Omega_F^+) \
\stackrel{\overline{\pi_\Omega}}{\mathop{\longrightarrow}\limits_{^{^\sim}}}
\ {\mathcal C}_H[\Omega _F^+]
\end{equation}
of (highest weight) categories.
\end{thm}
 This means, given
$\omega\in\Omega^+_F$ and $M\in\sC_G[\Lambda^+_F]$, there are
isomorphisms
\begin{equation}\label{Extisos}\begin{aligned}
\Ext^\bullet_H(\Delta_H(\omega),\pi_\Omega M)
&\cong\Ext^\bullet_{\sC_H[\Omega^+_F]}(\Delta_H(\omega),\overline{\pi_\Omega}j^*M)\\
&\cong\Ext^\bullet_{\sC_G(\Omega^+_F)}(j^*\Delta(\omega),j^*M)\\
&\cong\Ext^\bullet_{\sC_G[\Gamma^+_F]}(\Delta(\omega),M)\\
&\cong\Ext^\bullet_G(\Delta(\omega),M).
\end{aligned}\end{equation}
The last isomorphism follows since, at the derived category level,
$\sC_G[\Gamma^+_F]$ fully embeds into the category of rational
$G$-modules.  The same argument applies for the first isomorphism.
See \cite[Thm. 3.9]{CPS1}. The second isomorphism is a consequence
of Thm. \ref{maintheorem} and the fact, noted above, that
$j^*\Delta(\omega)$ is the standard object in the highest weight
category ${\mathcal C}_G(\Omega^+_F)$. The third isomorphism
follows from \cite[Lemma 6]{CPS6}; it is a general fact involving
a quotient functor $j^*:\sC_G[\Gamma^+]\to\sc(\Omega^+)$ for any
arbitrary coideal in a finite weight poset $\Gamma^+$.

Similarly, there are isomorphisms
\begin{equation}\label{Extisos2}\begin{aligned}
\Ext^\bullet_H(\pi_\Omega M,\nabla_H(\omega))&\cong
\Ext^\bullet_{\sC_H[\Omega^+_F]}(\overline{\pi_\Omega}j^*M,
\nabla_H(\omega))\\
&\cong\Ext^\bullet_{\sC_G(\Omega^+_F)}(j^*M,j^*\nabla(\omega))\\
&\cong\Ext^\bullet_{\sC_G[\Gamma^+_F]}(M,\nabla(\omega))\\
&\cong\Ext^\bullet_G(M,\nabla(\omega)).
\end{aligned}\end{equation}

In particular, in both (\ref{Extisos}) and (\ref{Extisos2}), we
can take $M=L(\tau)$, for $\tau\in \Omega^+_F$, to obtain
\begin{equation}\label{transfer1}
\Ext^\bullet_G(\Delta(\omega),L(\tau))\cong\Ext^\bullet_H(\Delta_H(\omega),
L_H(\tau))
\end{equation}
 and
\begin{equation}\label{transfer2}
\Ext^\bullet_G(L(\tau),\nabla(\omega))\cong\Ext^\bullet_H(L(\tau),\nabla(\omega)).
\end{equation}
See \cite[pp. 224--226]{CPS6} for further discussion, including
the relation of these $\Ext$-groups with Kazhdan-Lusztig
polynomials. In \S\S\S\ref{hemmersection},\ref{example},
 \ref{extoneisomorhismandhigherext}, we will take up the
 issue of comparing $\Ext$-groups between {\it
irreducible} modules.

\section{Some elementary applications}\label{kunnethformula}
  As noted in \cite[\S6]{CPS6}, the proof of Thm. \ref{maintheorem}
relies entirely on properties of highest weight categories, and so
it remains valid in other contexts, such as quantum enveloping
algebras and $q$-Schur algebras at a root of unity.   In this
section, we present various illustrations. These results will not
be used elsewhere in this paper.

First, we introduce some standard notation (which will be
important also in \S\ref{example}). If
$\lambda=(\lambda_1,\cdots,\lambda_n)$ is a partition of $r$ we
write $\lambda\vdash r$; let $\Lambda^+(r)$ be the set of all
partitions of $r$, and $\Lambda^+(n,r)$ the set of partitions with
at most $n$ nonzero parts.  If $\lambda\vdash r$, write
$|\lambda|=r$. The set $\Lambda^+(r)$ is regarded as a poset,
using the dominance ordering $ \unlhd$: $\lambda\
\unlhd\mu\iff\sum_{i=1}^j\lambda_i\leq\sum_{i=1}^j\mu_i$ for all
$j$. Then $\Lambda^+(n,r)$ is a coideal in $\Lambda^+(r)$. For
$\lambda\in\Lambda^+(r)$, let $\lambda'$ denote the dual
partition; thus, $\lambda_i'=\#\{\lambda_j\geq i\}$.

Let ${\mathbb E}$ be a vector space over the field of real numbers
$\mathbb R$ of dimension $n$. Assume that $\mathbb E$ is a
Euclidean space with orthonormal basis $\epsilon_1,\cdots,
\epsilon_n$. We consider the root system $\Phi$ of type $A_{n-1}$
of $\mathbb E$ which has roots $\epsilon_i-\epsilon_j$, $1\leq
i\not=j\leq n$. We label the roots, etc. in the standard way,
e.~g., as in \cite[p. 250]{Bo}. Thus, $\Phi$ has a simple set
$\Pi=\{\alpha_1,\cdots,\alpha_{n-1}\}$ of roots, putting
$\alpha_i-\epsilon_i-\epsilon_{i+1}$, $1\leq i<n$. We can identify
the $\mathbb Z$-lattice generated by $\epsilon_1,\cdots,
\epsilon_n$ with the character group $X=X(T)$ (written additively)
of the $n$-dimensional torus $T={\mathbb G}_m^{\times n}$, setting
$\epsilon_i:T\to {\mathbb G}_m$ to be the projection onto the
$i$-th factor. With this notation, the ``fundamental dominant
weights" are defined by putting
$\varpi_i=\epsilon_1+\cdots+\epsilon_i$, $i=1,\cdots, n-1$. Using
the $\mathbb Z$-basis $\epsilon_1,\cdots,\epsilon_n$ for $X$, we
can identify $X$ as the abelian group of sequences $(\lambda_1,
\cdots,\lambda_n)$ of integers.  Then $X^+$ identifies with the
set of all $(\lambda_1,\cdots,\lambda_n)$ satisfying
$\lambda_1\geq\cdots\geq\lambda_n$.  In this way,
$\Lambda^+(n,r)\subset X^+$.

For $1\leq d<n$, $\Pi_d:=\Pi\backslash\{\alpha_d\}$ is a base for
a subroot system $\Phi_d$ of type $A_{j-1}\times A_{n-1-j}$. For
$\lambda,\mu\in\Lambda^+(n,r)$, write $\lambda=_d\mu$ provided
that $\lambda-\mu\in{\mathbb Z}\Phi_d$. Also, put
$\lambda^{[d]}=(\lambda_1,\cdots,\lambda_d)$.  The following
result is immediately verified.

\begin{lem}\label{cosetlemma} For $\lambda,\mu\in\Lambda^+(n,r)$, $\lambda=_d\mu$ if
and only if $|\lambda^{[d]}|=|\mu^{[d]}|$.\end{lem}

Now, using Thm. \ref{maintheorem}, we easily obtain expressions
for various $\Ext$-groups involving algebraic groups (or quantum
groups) of type $A$.  In turn, the answers can be given in terms
of Schur algebras $S(n,r)$ (or $q$-Schur algebras $S_q(n,r)$). We
illustrate these results in several cases, referring the reader to
\cite{CPS6} for an explanation of the (very standard) notation,
and further references. In the first four examples, we work with
$q$-Schur algebras, letting  $q\in k$ be a primitive $\ell$th root
of unity.

\medskip\noindent
{\bf Example 1: Row removal.} Suppose
$\lambda,\mu\in\Lambda^+(n,r)$ satisfy $\lambda_1=\mu_1$. Let
$\overline\lambda=(\lambda_2,\cdots,\lambda_n),\overline\mu=(\mu_2,
\cdots,\mu_n)\in\Lambda^+(n-1,r')$, where $r'=r-\lambda_1$. Then
for $(M,N)\in\{(\Delta,L),(L,\nabla),
(\Delta,\Delta),(\nabla,\nabla)\}$, we have
$$\Ext^\bullet_{S_q(n,r)}(M(\lambda),N(\lambda))\cong
\Ext^\bullet_{S_q(n-1,r')}(M(\overline\lambda),N(\overline\mu)).$$
In fact, this follows immediately from Thm. \ref{maintheorem},
using Lemma \ref{cosetlemma} with $d=1$. This result has already
been observed in \cite[\S\S4,6]{CPS6}.

\medskip\noindent
{\bf Example 2: Removal of several rows.} Keep the notation of
Example 1 above, but replace $\Phi_1$ by $\Phi_d$ for some $d$
satisfying $1<d<n$. Now given $\lambda,\mu\in\Lambda^+(n,r)$,
suppose that $|\lambda^{[d]}|=|\mu^{[d]}|=r'$ and set $r''=r-r'$.
With $(M,N)$ as before,
$$\begin{aligned}\Ext^\bullet_{S_q(n,r)}&(M(\lambda),N(\lambda))\cong\\
&\Ext^\bullet_{S_q(d,r')}(M(\lambda^{[d]}),N(\mu^{[d]})) \otimes
\Ext^\bullet_{S_q(n-d,r'')}(M(\lambda\backslash\lambda^{[d]}),
N(\mu\backslash\mu^{[d]})),\end{aligned}$$
 as graded vectors spaces.  We have written
 $\lambda\backslash\lambda^{[d]}$ for
 $(\lambda_{d+1},\cdots,\lambda_n)$ and similarly for $\mu$.
In checking these facts, it is useful to make use of the fact
that, for positive integers $a,b$,
$\Ext^\bullet_{S_q(a,b)}(M(\lambda),N(\lambda))\cong
\Ext^\bullet_{S_q(b,b)}(M(\lambda),N(\lambda))$ for
$\lambda,\mu\in\Lambda^+(a,b)$.  (On the right-hand side, we
interpret $\lambda,\mu\in\Lambda^+(b)$.) Suitably interpreted in
the context of $q$-Schur algebras, this isomorphism is just that
given in (\ref{transfer1}) or (\ref{transfer2}); see \cite[Thm.
6.1]{PS}.

\medskip\noindent
{\bf Example 3: Passage to Hecke algebras.} Define
     $e=\begin{cases} p\quad {\text{\rm if $\ell=1$}};\\
                  \ell\quad{\text{\rm
otherwise}},\end{cases}$
  where $k$ has characteristic $p$.
  Let $\sH_r$ be the Hecke algebra associated to the symmetric
  group ${\mathfrak S}_r$ and the parameter $q$.  We use the
  notation of \cite{PS}.
 For $\lambda\vdash r$, let
$S_\lambda$ denote the corresponding Specht module for $\sH_r$.
 Let $m$ be a non-negative integer satisfying $m\leq e-3$.
 Consider the setting of Example 2 with $n=r$, and let
 $\lambda,\mu\in\Lambda^+(r)$ satisfy $\lambda^{[d]}=\mu^{[d]}$.
 Then
$$\Ext^m_{\sH_r}(S_\lambda,S_\mu)=\bigoplus_{i=0}^m\Ext^i_{\sH_s}(S_{\lambda^{[d]}},
      S_{\mu^{[d]}})\otimes\Ext^{m-i}_{\sH_t}(S_{\lambda\backslash\lambda^{[d]}},
      S_{\mu\backslash\mu^{[d]}}).$$
This result is immediate from \cite[Thm. 4.6(ii)]{PS}.

\medskip\noindent
{\bf Example 4: Column removal.} Fix $d$, $1\leq d<n$, and suppose
that $\lambda,\mu\in\Lambda^+(n,r)$ satisfy $\lambda'=_d\mu'$.
Then we can use Example 2 to give a K\"unneth factorization of
$\Ext^\bullet_{S_q(n,r)}(\Delta(\lambda),L(\mu))$, after observing
that
\begin{equation}\label{tilt}
\Ext^\bullet_{S_q(r,r)}(\Delta(\lambda),\Delta(\mu))\cong\Ext^\bullet_{S_q(
r,r)}(\Delta(\mu'), \Delta(\lambda')).
\end{equation}
To see this, we use some tilting module results, summarized (with
references) in \cite[\S3]{CPS6}. Let $Y$ be a complete tilting
module for $S_q(r,r)$--mod. Set $E_q(r,r)=\End_{S_q(r,r)}(Y)$.
Form the contravariant tilting functor
$T:S_q(r,r){\text{\rm--mod}}\to E_q(r,r){\text{\rm --mod}}$ given
by $T(M):=\Hom_{S_q(r,r)}(M,Y)$. The category $E_q(r,r){\text{\rm
--mod}}$ is equivalent to the category $S_q(r,r)$--mod, so we
identify these two categories. (Caution: the equivalence is {\it
not} by means of the tilting functor $T$.) With this
identification,  $T\Delta(\lambda)\cong \Delta(\lambda')$. If
$P_\bullet\to \Delta(\lambda)\to 0$ is a projective resolution,
then $0\to\Delta(\lambda')\to TP_\bullet$ is a resolution of
$\Delta(\lambda')\cong T\Delta(\lambda)$ by tilting modules; see
\cite[Lemmas 1.4\&1.5]{PS}. Since tilting modules are acyclic for
the functor $\Hom_{E_q(r,r)}(\Delta(\mu'),-)$, (\ref{tilt})
follows from the isomorphism
$$\Hom_{E_q(r,r)}(\Delta(\mu'),TP_\bullet)\cong
\Hom_{S_q(r,r)}(P_\bullet,\Hom_{E_q(r,r)}(\Delta(\mu'),Y))$$
   (see \cite[Prop.
1.2]{PS}) and the fact that $\Hom_{E_q(r,r)}(\Delta(\mu'),Y)\cong
\Delta(\mu)$. More generally, this argument shows that $T$ defines
a contravariant equivalence from the exact subcategory of
$S_q(r,r)$-modules with a $\Delta$-filtration to the similar
category for $E_q(r,r)$--mod.

 \medskip
In Example 2, the isomorphism
$$\Ext^\bullet_{S_q(n,r)}(\nabla(\lambda),\nabla(\mu))\cong
\Ext^\bullet_{S_q(d,r')}(\nabla_{S_q(d,r')}(\lambda),\nabla_{S_q(d,r')}(\mu))$$
is already essentially proved in \cite[Formula (17), p. 91]{D2},
without the explicit interpretation in terms of the decomposition
as a tensor product.

 The K\"unneth results above
  have only been worked out in type $A$. More generally, {\it in
other types} in the case of a reductive algebraic (or quantum)
group $G$, a similar such tensor factorization would exist
whenever the semisimple part of the Levi subgroup $H$ decomposes
into a product of smaller reductive groups. However, it is not
likely that the  elegant interpretation   in type $A$ given above
holds in all other cases. Our last example below indicates that
this does occur in type $C_n$ at least.

\medskip\noindent
 {\bf  Example 5: Type $C$.}  Suppose that $G\cong Sp_{2n}(k)$ has
type $C_n$.  We use the notation of \cite[p.254-255]{Bo}
(replacing $l$ there by $n$) as far as listing the set $\Pi$ of
simple roots, etc. Thus, $X=X(T)\cong{\mathbb
Z}\epsilon_1\oplus\cdots\oplus{\mathbb Z}\epsilon_n$, and the
fundamental dominant weights $\varpi_i$, $1\leq i\leq n$, are
given by $\varpi=\epsilon_1+\cdots+\epsilon_i$. The set
$\Lambda^+(n,\bullet)$ of all partitions $\lambda$ of length at
most $n$ indexes the dominant weights $X^+$. Just as for type $A$,
for $1\leq d\leq n$, define $\Phi_d$ as the subroot system of type
$A_{d-1}\times C_{n-d-1}$ of $\Phi$ with simple roots
$\Pi_d:=\Pi\backslash\{\alpha_d\}$. Put $\lambda=_d\mu$ if and
only if $\lambda-\mu\in{\mathbb Z}\Phi_d$.  Then, by analogy with
Lemma \ref{cosetlemma}, given
$\lambda,\mu\in\Lambda^+(n,\bullet)$, we verify that
$\lambda=_d\mu$ if and only if $|\lambda^{[d]}|=|\mu^{[d]}|$ and
$|\lambda|\equiv|\mu|$ mod$\,2$. If $H$ is the Levi subgroup of
$G$ corresponding to $\Phi_d$, then $H\cong GL_d(k)\times
Sp_{2(n-1-d)}(k)$.  Thus, if $\lambda=_d\mu$, the groups
$\Ext^\bullet_G(M(\lambda),N(\lambda))$ are isomorphic as a graded
vector space to
$$\Ext^\bullet_{GL_d(k)}(M(\lambda^{[d]}),N(\mu^{[d]}))\otimes
\Ext^\bullet_{Sp_{2g}(k)}(M(\lambda\backslash\lambda^{[d]}),N(\mu\backslash\mu^{[d]})),$$
where $g=n-d-1$ and $(M,N)\in\{(\Delta,L),(L,\nabla),
(\Delta,\Delta),(\nabla,\nabla)\}$.  We leave further details to
the interested reader.

\medskip\medskip
\begin{center}{\large\bf PART II: Character formulas}\end{center}
\addcontentsline{toc}{section}{\bf Part II: Character formulas}
\medskip

We study various conditions equivalent to the validity of the
Lusztig character formula.  We arrange our formulations, in this
section and the next, so that we can discuss cases where the
formula might hold, even though the full Lusztig conjecture might
not hold, or might not be known to hold.

 \section{The Lusztig character formula}\label{charformulaholds}

  Let $G$ be a
simple and simply connected algebraic group over $k$. (The results
below and in Part III easily extend to the case of a general
reductive group.) In addition, let $\rho\in X^+$ be defined by
$\langle\rho,\alpha^\vee\rangle=1$ for all simple roots
$\alpha\in\Pi$.  Let $W=N(T)/T$ be the Weyl group of $G$ and let
$W_p:=W\\ltimes p{\mathbb Z}\Phi$ be the affine Weyl group,
generated by $W$ and the normal subgroup $p{\mathbb Z}\Phi$ of
translations by $p$-multiples of roots.  It acts naturally (and
faithfully) on the space ${\mathbb E}={\mathbb R}\otimes X$. For
$\alpha\in\Phi$, $n\in{\mathbb Z}$, $s_{\alpha,np}\in W_p$ is
defined as an operator on ${\mathbb E}$ by $s_{\alpha,np}(u)=u -
(\langle u,\alpha^\vee\rangle -np)\alpha$, $u\in {\mathbb E}$. Set
$\wSp=\{s_\alpha\,|\,\alpha\in\Pi\}\cup\{s_{\alpha_o,-p}\}\subset
W_p$, where  $\alpha_o\in\Phi$ is the maximal short root. Then
$(W_p,\wSp)$ is a Coxeter system; let  $l:W_p\to{\mathbb Z}$ be the
length function of $W_p$ defined by $\wSp$. In place of the natural
action of $W_p$ on ${\mathbb E}$, we prefer the ``dot" action given
by $w\cdot u= w(u+\rho)-\rho$, $w\in W_p$, $u\in{\mathbb E}$. Given
$\alpha\in\Phi$, $n\in{\mathbb Z}$, let $H_{\alpha,np}
=\{u\in{\mathbb E}\,|\,\langle u +\rho,\alpha^\vee\rangle=pn\}.$
Thus, $s_{\alpha,np}\in W_p$ is a reflection about the hyperplane
$H_{\alpha, np}$. Any connected component $C$ of ${\mathbb
E}\backslash \bigcup_{\alpha\in\Phi, n\in{\mathbb Z}}H_{\alpha,np}$
is an alcove of $W_p$.  The closure $\overline C$ a fundamental
domain for the dot action of $W_p$ on $\mathbb E$.

A weight $\lambda\in X$ is called {\it regular} provided that
$\lambda$ belongs to some alcove $C$ of $W_p$, i.~e., provided
$\langle\lambda +\rho,\alpha^\vee\rangle\not\equiv 0$ mod $p$ for
all $\alpha\in\Phi$. If $h$ denotes the Coxeter number of $G$,
  regular weights exist if and only if $p\geq h$. In
particular, if $p\geq h$, $-2\rho=w_o\cdot 0$ is a regular weight,
where $w_o$ denotes the long word in $W$.  If $\Sigma\subset
{\mathbb E}$, let $\Sigma_{\text{\rm reg}}$ be the set of all
regular elements in $\Sigma$.

For the rest of this section, assume that $p\geq h$. By
definition, {\it the Lusztig character formula holds for}
$\lambda\in X^+$ (in the category ${\mathcal C}_G$ of finite
dimensional rational $G$-modules) provided
$\lambda=x\cdot(-2\rho)$ for some $x\in W_p$ and the irreducible
$G$-module $L(\lambda)$ has formal character given by
\begin{equation}\label{lusztigcharacterformula}
\text{\rm ch}\,L(\lambda)=\sum_{y\cdot({-2\rho})\in
X^+}(-1)^{\l(x)-\l(y)}P_{y,x}(1)\text{\rm
ch}\,\Delta(y\cdot(-2\rho)).
\end{equation}
In this formula, $P_{y,x}$ is the Kazhdan-Lusztig
polynomial\footnote{The polynomials $P_{y,x}$ are associated to
the generic Hecke algebra $\mathcal H$ for the Coxeter system
$(W_p,\wSp)$. Each $P_{y,x}$ is, in fact, a polynomial in $u:=t^2$
and $P_{x,y}(1)$ means the specialization at $u=1$ (i.~e., at
$t=-1$). The introduction of the variable $t=\sqrt{u}$ is only
important in this paper in the proof of Thm.
\ref{lusztigequiv}(b),(c) below, involving the theory of enriched
Grothendieck groups.} associated to the pair $(y,x)\in W_p\times
W_p$. In particular, $P_{y,x}=0$ unless $y\leq w$ in the
Bruhat-Chevalley ordering on $W_p$. For any $\xi\in X^+$,
$$\text{\rm ch}\,\Delta(\xi)=\frac{\sum_{w\in
W}e^{w\cdot\xi}}{\sum_{w\in W}e^{w\cdot 0}}$$
  by Weyl's character
formula.
 Given a subset $\Sigma\subseteq W_p\cdot 0=W_p\cdot(-2\rho)$, we
say that the Lusztig character formula holds for $\Sigma$ it if
holds for every $\lambda\in\Sigma$.

The famous Lusztig conjecture for modular characters \cite{L}
states that the Lusztig character formula holds for $W_p\cdot
0\cap \Jan$, where
\begin{equation}\label{jantzenregion} \Jan=\{\lambda\in
X^+\,|\,\langle\lambda+\rho,\alpha^\vee\rangle\leq
p(p-h+2),\quad\forall \alpha\in\Phi\}
\end{equation}
is the Jantzen region. Of course, to test if $\lambda\in\Jan$, it
suffices to check the inequality in (\ref{jantzenregion}) just for
$\alpha=\alpha_o$. For a group $G$ with a fixed root system
$\Phi$, the conjecture holds provided $p$ is large enough,
through, in general, no sufficient bound on $p$ is known. See
\cite{AJS} and, for a general discussion, \cite[\S8]{T}. The
Lusztig conjecture holds generally in types $A_n$, $n\leq 3$, and
$B_2$, $G_2$ \cite{Jan}. In addition, it is known to hold in type
$A_4$ for $p=5, 7$; see the discussion in \cite{S2}.

However, independently of the validity of the Lusztig conjecture,
ideals $\Gamma$ for which the Lusztig character formula holds for
$W_p\cdot 0\cap\Gamma$ occur naturally, as Thm. \ref{lusztigequiv}
below suggests. A simple (almost trivial) example would be to take
$\Gamma=\{\lambda\in X^+\,|\,\lambda\leq s_{\alpha_o,p}\cdot 0\}$.

Let  $C^+\subset{\mathbb E}$ be the alcove containing $0$, and put
$C^-=w_o\cdot C^+$ (the alcove containing $-2\rho$).  Let
$\lambda\in X\cap C^-$ and suppose the Lusztig character formula
holds for $x\cdot(-2\rho)$, $x\in W_p$. Using translation,
\begin{equation}\label{genlusztigformula}
\text{\rm ch}\,L(x\cdot\lambda)=\sum_{y\cdot\lambda\in
X^+}(-1)^{\l(x)-\l(y)}P_{y,x}(1)\text{\rm
ch}\,\Delta(y\cdot\lambda).\end{equation}
  Thus, we can say that the Lusztig character formula holds for
  $x\cdot\lambda$. Suppose that $\Sigma\subset
X^+$ has the property that $\Sigma_{\text{\rm reg}}$ is a union of
sets of the form $C\cap X^+$, $C$ an alcove. (For example, $\Jan$
has this property.) Then, if the
  Lusztig character formula holds for $W_p\cdot 0\cap\Sigma$,
  it also holds for $\Sigma_{\text{\rm reg}}$.\footnote{Suppose
  $\lambda'\in X\cap \overline{C^-}$ is not
  regular. Let ${\mathcal D}_{\lambda'}$ be set of distinguished left coset representatives
  in $W_p$ of the stabilizer in $W_p$ of $\lambda'$ under the dot action.
  Suppose that (\ref{genlusztigformula}) holds  for a given
  $x\in{\mathcal D}_\lambda$. Then, using Jantzen translation,
  $\ch L(\x\cdot\lambda')=\sum_{y\cdot\lambda'\in
  X^+ }(-1)^{\l(x)-\l(y)}P_{y,x}(1)\text{\rm
ch}\,\Delta(y\cdot\lambda')$.  In this way, it would be possible
to speak of the Lusztig character formula holding for an ideal
$\Gamma$, not just $\Gamma_{\text{\rm reg}}$, though we will not
use this more general terminology.}

Let $X^+_1$ be the set of $p$-restricted dominant weights. A
conjecture of Kato \cite[Conj. 5.5]{K} asserts that the Lusztig
character formula holds for the set $X_1^+\cap W_p\cdot 0$ if
$p\geq h$. Observe that $X_1^+\subseteq \Jan$ if and only if
$p\geq 2h-3$. If his conjecture holds, then (the argument in)
\cite[(5.4)]{K} also establishes that
(\ref{lusztigcharacterformula}) holds for any
$L(\lambda)=L(\tau)\otimes L(\xi)^{(1)}$, where $\tau\in
X^+_{1\,\text{\rm reg}}$ and  $\xi\in X^+$ has the property that
$\Delta(\xi)\cong L(\xi)$.\footnote{In general, such weights
$\tau+p\xi$ need not always lie in $\Jan$.} Here, given a rational
$G$-module $V$, $V^{(1)}$ denotes the ``twist" of $V$ through the
Frobenius morphism ${\text{Fr}}:G\to G$.

The following result will be useful in \S\ref{example}. Although
the statement of (a) is purely combinatorial, the authors are
unaware of a purely combinatorial proof. It would also be
interesting to remove the various restrictions in the statement
this result.

\begin{lem}\label{amazinglemma}Assume that $p\geq h$ and that the
root system of $G$ is simply laced. In addition, if $G$ has type
$E_6, E_7, E_8$ assume that $p>h+1$.  Suppose that $x\in W_p$ and
$x\cdot(-2\rho)=p\lambda$ for $\lambda\in X^+$.

(a) We have $$ \text{\rm ch}\,\Delta(\lambda)^{(1)}=
\sum_{y\cdot({-2\rho})\in
X^+}(-1)^{\l(x)-\l(y)}P_{y,x}(1)\text{\rm
ch}\,\Delta(y\cdot(-2\rho)).
$$

(b) If $\Delta(\lambda)\not=L(\lambda)$, the Lusztig character
formula (\ref{lusztigcharacterformula}) does not hold for
$p\lambda=x\cdot(-2\rho)$.
\end{lem}
\begin{proof} Let $U_q$ be the quantum enveloping algebra (integral
form) at $q=\sqrt[\leftroot{2} \uproot{2} p]{1}$ of the same type
as $G$. Let ${\mathcal C}_q$ be the category of finite dimensional
integral type 1 modules for $U_q$. It is a highest weight category
with irreducible (resp. standard) modules $L_q(\xi)$ (resp.,
$\Delta_q(\xi)$), $\xi\in X^+$. Furthermore, for any
$\xi=x\cdot(-2\rho)\in X^+$, the Lusztig character formula holds
for $\xi$ in the category ${\mathcal C}_q$; in other words,
(\ref{lusztigcharacterformula}) remains valid if each ``$L$"
(resp., ``$\Delta$") is replaced by ``$L_q$" (resp.,
``$\Delta_q$") This last assertion follows essentially from the
validity of Lusztig's quantum group conjecture: see \cite[\S7]{T}
for a detailed discussion and further references. It is here that
the restriction that $p>h+1$ is required in the exceptional
cases---see the (equivalent) case--by--case formulation in
\cite[Thm. 7.1]{T}. Of course, $\text{\rm
ch}\,\Delta_q(\xi)=\text{\rm ch}\,\Delta(\xi)$ for any $\xi\in
X^+$.

Let $\sigma:U_q\to U$ be the Frobenius morphism from $U_q$ to the
universal enveloping algebra over ${\mathbb C}$ of the complex Lie
algebra $\mathfrak g$ having the same type as $G$. Then
$\Delta_q(p\lambda)=\sigma^*L_{\mathbb C}(\lambda)$, the pull-back
through $\sigma$ of the complex irreducible $\mathfrak g$-module
of highest weight $\lambda$. Hence,
$\ch\,\Delta_q(p\lambda)=\ch\,\Delta(\lambda)^{(1)}$, and (a)
follows.

Since $L(p\lambda)\cong L(\lambda)^{(1)}$, (b) also follows.
\end{proof}

\begin{rem}See \cite{T} for a discussion of the status of the Lusztig
quantum group conjecture used in the above proof in the non-simply
laced case. Published results giving a character formula require
$p\geq r_0$, where $r_0$ is sufficiently large, depending on the
root system, but unknown; see \cite[Thm. 7.1]{T} and remarks
following it. However, \cite{T} also notes an improvement to $p>h$
in all cases using unpublished work of Soergel.  A similar
improvement in all cases might be based on the category
equivalences given in \cite{ABG}, though specific details leading
from the final result in \cite{ABG} to a character formula are not
included. In our paper, the above lemma is only used in type $A$.
In fact, in this case, as noted in \cite{T}, the equivalence
between representations of the quantum enveloping algebra at an
$\ell$th root of unity and the corresponding representations for
affine Lie algebras holds for all $\ell$, without any $\ell\geq h$
restriction.

\end{rem}
 \section{Equivalent
 conditions}\label{equivalentconditionssection}

Assume $G$ is simple (any type), simply connected and $p\geq h$.
The faces of an alcove $C$ are naturally labelled by elements
 of $\wSp$: if $C=C^-$, the faces are already labelled by $\wSp$ by definition; otherwise,
 a face $F$ is $W_p$ dot-conjugate to a unique face of $C^-$, and we assign to
 $F$ the corresponding element $s\in \wSp$. Thus, for $s\in \wSp$, we speak of an $s$-face of $C$.
 For
 $\lambda,\mu\in X^+_{\text{\rm reg}}$ and $s\in \wSp$, write $\mu=\lambda s$
 provided: (1)
  $\lambda,\mu$ lie in adjacent alcoves $C,C'$ separated by an $s$-face, and (2)
  $\lambda$ is the reflection of $\mu$ through that $s$-face.
  In other words, $\lambda=s_{\beta,np}\cdot\mu$, where the
  common $s$-face of $C,C'$ lies in the hyperplane $H_{\beta,np}$.
Then $\lambda$ and $\lambda s$ are called {\it adjacent}. If
$\lambda=w\cdot(-2\rho)$, $w\in W_p$, then $\lambda
s=ws\cdot(-2\rho)$. Also, $\lambda s>\lambda\iff l(ws)=l(w)+1$,
and $\lambda=w\cdot(-2\rho)\in X^+\iff w=w_oy$, with $y$ a
``distinguished" right coset representative of $W$ in $W_p$,
i.~e., $l(w)=l(w_o)+l(y)$. In particular, if $y$ is distinguished
and $y=s_1\cdots s_n$, $s_i\in \wSp$, then each product $s_1\cdots
s_m$, $m\leq n$, is also a distinguished distinguished right coset
representative. Thus,
$$0=w_o\cdot(-2\rho)<w_os_1\cdot(-2\rho)<\cdots< w_os_1\cdots
s_n\cdot(-2\rho)$$ provides a ``path" of adjacent dominant weights
from $0$ to $\lambda$.

In the following  theorem, part (a)  plays a basic role in the
sequel.  Part (b) gives a way to construct ideals $\Gamma$ for
which the Lusztig character formula holds.  However, it and part
(c)  will not be used in the paper,\footnote{More precisely, (c)
makes a brief cameo appearance in \S\ref{example}.} so that the
proofs are rather brief, and rely heavily on some of the machinery
developed in \cite{CPS2}, \cite{CPS4}.

\begin{thm} \label{lusztigequiv}  Let $\Gamma$ be a finite ideal in the poset $(X^+,\leq)$.

 (a) The Lusztig character formula holds
for $\Gamma\cap W_p\cdot 0$ if and only if
\begin{equation}\label{Extcondition} \Ext^1_G(L(\lambda),L(\lambda
s))\not=0
\end{equation}
whenever $\lambda<\lambda s$ both belong to $\Gamma\cap W_p\cdot
0$.

(b) Suppose that there exists $\tau\in X^+\cap W_p\cdot 0$ such
that $\Gamma=\{\xi\in X^+\,|\, \xi\leq\tau\,\}$ and form the ideal
$\widetilde\Gamma=\{\xi\in\Gamma\,|\,\xi<\tau\}$ of $\Gamma$.
 Suppose that $s\in \wSp$ and $\tau s\in\widetilde\Gamma$. Assume that the Lusztig
character formula holds for $\widetilde\Gamma\cap W_p\cdot 0$.
Then the Lusztig character formula holds for $\Gamma\cap W_p\cdot
0$ if and only if $\Ext^1_G(L(\tau s),L(\tau))\not=0$.\footnote{In
the statement of (b) and (c), we could replace $\tau$ and $\xi$ be
regular weights, using a standard translation functor argument.}

(c) Given $\xi\in X^+\cap W_p\cdot 0$, the Lusztig character
formula holds for $\xi$ provided there is a path
$0=\xi_0<\xi_1<\cdots<\xi_m=\xi$ in $X^+\cap W_p\cdot 0$ with
$\xi_{_i-1}$ adjacent to $\xi_i$ for $0<i\leq m$ and with each
$\Ext^1_G(L(\xi_{i-1}),L(\xi_i))\not=0$.
\end{thm}
\begin{proof} Part (a) is proved in \cite[Thm. 5.3]{CPS2} for
$\Gamma=\Jan$, but
 the argument holds equally well for any ideal $\Gamma$ in $X^+$, making use of the
 discussion in the paragraph immediately above the statement of the theorem.

Although (c) implies (b), we prove (b) first, then augment the
argument to obtain (c). In in the proof of (b), we will need to
use the theory (and notation) of pre-Hecke operators discussed in
\cite[\S5]{CPS2} and \cite[\S4]{CPS4} (the latter reference
contains many more details). Let ${\mathcal C}_{G,0}$ be the
category of all finite dimensional rational $G$-modules with
composition factors $L(\lambda)$ for some $\lambda\in W_p\cdot 0$.
Then ${\mathcal C}_{G,0}$ is a highest weight category with poset
$(W_p\cdot 0,\uparrow)$, letting $\uparrow$ be the Jantzen ``up
arrow" poset structure. Form the exact, additive functor
$\Theta_s:{\mathcal C}_{G,0}\to{\mathcal C}_{G,0}$ (the composite
of two Janzten translation functors) \cite[p. 88]{CPS4}.  If $M$
has the form $\Delta(\lambda), \nabla(\lambda)$, or $L(\lambda)$
for some $\lambda\in W_p\cdot 0$, then let $\beta_sM$ equal to the
complex
 $$ 0\to M\to \beta_sM\to M\to 0,$$
 defined in \cite[(4.8.2)]{CPS4}.  It is concentrated in degrees
 $-1$, 0, and $+1$.  When $\lambda < \lambda s$, $\beta_sL(\lambda)$ is isomorphic in
 the bounded derived category
 $D^b({\mathcal C}_{G,0})$ to its cohomology, which is concentrated in
 degree 0. (We will take $\lambda s = \tau$ in the notation of (b).)
  From the discussion in \cite[pp. 89--90]{CPS4}, if
 $\lambda \in\widetilde \Gamma$, $\beta_sL(\lambda)\in {\mathcal
 E}^L\cap{\mathcal E}^R$. See also \cite[Remark 2.2]{CPS4}.
Here ${\mathcal E}^L\cap{\mathcal E}^R$ is a subset of objects in
$D^b({\mathcal C}_{G,0})$ defined by an even-odd vanishing
property similar to that described for irreducible modules in
(\ref{KLtheory}) below.
   The condition $\Ext^1_G(L(\tau
 s),L(\tau))\not=0$ implies that $L(\tau)$ is a direct summand of
 $\beta_sL(\tau s)$; see the end of the proof of ``$(c)\implies (b)$" of
  \cite[Thm. 5.5]{CPS4}.  Therefore, $L(\tau)\in{\mathcal
 E}^L\cap{\mathcal E}^R$.  This implies that ${\mathcal
 C}_0[\Gamma]$ has a Kazhdan-Lusztig theory (again, see
 (\ref{KLtheory}) below), and so the Lusztig character formula
 holds for $\Gamma$, as required.

The proof of (c) requires the use of the ``enriched Grothendieck
groups'' $K^L_0:=K^L_0({\mathcal C}_{G,0},l)$ and
$K^R_0:=K^R_0({\mathcal C}_{G,0},l)$ introduced in
\cite[\S2]{CPS2} (see also \cite[\S2]{CPS4}). In fact, these
groups can be associated to any highest weight category ${\mathcal
C}$, with weight poset $\Lambda$ and length function
$l:\Lambda\to{\mathbb Z}$.  As discussed in \cite[Prop.
2.3]{CPS2}, $K^L_0$ and $K^R_0$ are related by a non-degenerate
sesquilinear pairing $K^L_0\times K^R_0\to{\mathbb Z}[u,u^{-1}]$
(Laurent polynomials). When ${\mathcal C}_{G,0}[\Gamma\cap
W_p\cdot 0]$ has a Kazhdan-Lusztig theory \cite[Defn. 3.3]{CPS2}
(see also (\ref{KLtheory}) below), the group $K^L_0$ (resp.,
$K^R_0$) keeps track of the dimensions of the groups
$\Ext^\bullet_{{\mathcal C}_{G,0}}(L(\lambda),\nabla(\mu))$
(resp., $\Ext^\bullet_{{\mathcal
C}_{G,0}}(\Delta(\mu),L(\lambda))$) for $\lambda,\mu\in\Gamma\cap
W_p\cdot0$.

Let ${\mathcal H}$ denote the generic Hecke algebra of the Coxeter
system $(W_p,\wSp)$ over ${\mathbb Z}[t,t^{-1}]$. Then $\mathcal
H$ acts on the groups $K^L_0$ and $K^R_0$, using the $\beta_s$,
$s\in \wSp$, mentioned above \cite[Prop. 5.6]{CPS2}. The argument
for (b) shows, recursively, that all the $L(\xi_i)$'s may be
represented in both $K^L_0$ and $K^R_0$. Also, the discussion in
\cite[\S4]{CPS4} shows the resulting representations are
completely determined by Kazhdan-Lusztig polynomials.  Using the
``Euler characteristic map" $K^0_L\to{K}_0({\mathcal C}_{G,0})$
(the Grothendieck group of $G$--mod) obtained by specializing
$u\mapsto -1$ (and hence $t\mapsto 1$; cf. footnote 4) yields
(\ref{lusztigcharacterformula}) for $L(\xi)$.
\end{proof}

Let $\lambda\in X^+_{\text{\rm reg}}$ and suppose that
$\lambda<\lambda s\in X^+$ for some $s\in \wSp$. By \cite[II,
Prop. 7.21]{Jan}, $\dim\,\Ext^1_G(\Delta(\lambda),\Delta(\lambda
s)=1$. In simple terms, this means there is a non-split short
exact sequence
\begin{equation}\label{shortexactsequence1}
0\to \Delta(\lambda s)\to E\to \Delta(\lambda)\to 0
\end{equation}
of rational $G$-modules, and that, given any other non-split short
exact sequence $0\to\Delta(\lambda s)\to E'\to \Delta(\lambda)\to
0$, it is ``scalar equivalent" to (\ref{shortexactsequence1}) in
the sense that there is a commutative diagram
 \begin{equation}\label{deltacommdiagram}
\begin{array}{cccccccccc}
  0 & \begin{picture}(23,0)
\put(0,5){\vector(1,0){23}}
\end{picture}
  &\Delta(\lambda s)&
\begin{picture}(23,0)
\put(0,5){\vector(1,0){23}}
\end{picture}
   & E &
\begin{picture}(23,0)
\put(0,5){\vector(1,0){23}}
\end{picture}
  &\Delta(\lambda)&
\begin{picture}(23,0)
\put(0,5){\vector(1,0){23}}
\end{picture}
& 0
\\ &&
\begin{picture}(0,23)
\put(-17,7){$\scriptstyle{\alpha \, {\rm id}}$}
\put(0,21){\vector(0,-1){23}}
\end{picture}
  &&
\begin{picture}(0,23)
\put(0,21){\vector(0,-1){23}}
\end{picture}
&&
\begin{picture}(0,23)
\put(4,7){$\scriptstyle{{\rm id}}$} \put(0,21){\vector(0,-1){23}}
\end{picture} &&
\\
  0&
\begin{picture}(23,0)
\put(0,5){\vector(1,0){23}}
\end{picture}
&\Delta(\lambda s)&
\begin{picture}(23,0)
\put(0,5){\vector(1,0){23}}
\end{picture} &
E'&
\begin{picture}(23,0)
\put(0,5){\vector(1,0){23}}
\end{picture} &\Delta(\lambda)&
\begin{picture}(23,0)
\put(0,5){\vector(1,0){23}}
\end{picture} &0
\end{array}
\end{equation}
in which $0\not=\alpha\in k$.
 In general,
$$\dim\Ext^1_G(L(\lambda),L(\lambda s))\leq 1.$$
 If a
non-trivial extension of $L(\lambda)$ by $L(\lambda s)$ exists, it
must come, up to scalar equivalence, from
(\ref{shortexactsequence1})  by factoring out $\text{\rm
rad}(\Delta(\lambda s))$.  That is, the resulting extension of
$\Delta(\lambda)$ by $L(\lambda s)$ must split on $\text{\rm
rad}(\Delta(\lambda))$.  If it does, we get a non-split extension
of $L(\lambda)$ by $L(\lambda s)$; if it does not, there is no
such nontrivial extension. Observe that
\begin{equation}\label{transferofextone}
 \Ext^1_G(\Delta(\lambda),\Delta(\lambda
s))\cong\Ext^1_G(\Delta(\lambda),L(\lambda s)).
\end{equation}These statements can be easily
verified, using the long exact sequence of $\Ext^\bullet$. Note
that $\Ext^n_G(\Delta(\lambda),L(\nu))=0$ unless $\nu>\lambda$.

For each $\mu\in\Gamma\cap W_p\cdot 0$, let $P_\mu$ be a fixed
projective in the category ${\mathcal C}_G[\Gamma]$  such that the
projective cover $P(\mu)$ of $L(\mu)$ is a direct summand of
$P_\mu$, and all other indecomposable summands $P(\nu)$ satisfy
$\mu<\nu\in\Gamma$.\footnote{Such projective modules $P_\mu$ arise
naturally in many contexts, where $P(\mu)$ itself may not be
explicitly known. We allow summands $P(\nu)$ for $\nu\not\in
W_p\cdot 0$, though they will be irrelevant, since then $L(\nu)$
and $L(\mu)$  lie in different blocks. \par We only use that the
other summands $P(\nu)$ of $P_\mu$ satisfy  $\mu\leq\nu$ in this
paper, though the natural examples have $\mu<\nu$. That is, we
assume that $P(\mu)$ is a summand of $P_\mu$ with multiplicity
one.} Explicitly, $P_\mu\in {\mathcal C}_G[\Gamma]$, while
$P(\mu)\in {\mathcal C}_G[\Gamma\cap W_p\cdot 0]$, the category of
all finite dimensional rational $G$-modules with composition
factors of the form $L(\nu)$, $\nu\in\Gamma\cap W_p\cdot 0$. For
$\gamma\in\Gamma$, let $P_{\mu,<\gamma}$ be the largest quotient
module of $P_\mu$ with all composition factors $L(\nu)$,
$\nu<\gamma$.

The module $P_{\mu,<\gamma}$ has a filtration by submodules with
successive sections standard modules $\Delta(\tau)$, $\tau<\gamma$.

\begin{thm}\label{firstreduction}Let $\Gamma$ be a finite ideal
as above. For any given $\lambda,\lambda s\in\Gamma\cap W_p\cdot
0$ with $\lambda<\lambda s$, $s\in \wSp$, condition
(\ref{Extcondition}) holds if and only if, for each $\mu<\lambda$
in $\Gamma\cap W_p\cdot 0$, the following condition
 holds:
\begin{equation}\label{theconditionstar}
\left\{\parbox{280pt}{
 For any map $P_{\mu,<\lambda s}\to\Delta(\lambda)$ which
is not a split surjection, the induced map
$\Ext^1_G(\Delta(\lambda),\Delta(\lambda
s))\to\Ext^1_G(P_{\mu,<\lambda s},\Delta(\lambda s))$ is
zero.}\right.
\end{equation}
Thus, the Lusztig character formula holds for $\Gamma\cap W_p\cdot
0$ if and only if (\ref{theconditionstar}) holds for all
$\mu,\lambda,\lambda s\in\Gamma\cap W_p\cdot 0$.
\end{thm}

\begin{proof}We will make use of elementary recollement results; these are described
in \cite[\S2]{CPS6}, for example. Let $\Gamma'$ be an ideal in
$\Gamma$. If $M\in{\mathcal C}_G[\Gamma]$, $M_{\Gamma'}$ denotes
the largest quotient module of $M$ all of whose composition
factors lie in $\Gamma'$. If $P\in{\mathcal C}_G[\Gamma]$ is
projective, then $P_{\Gamma'}\in{\mathcal C}_G[\Gamma']$ is
projective. Moreover, the natural inclusion ${\mathcal
C}_G[\Gamma']\subseteq {\mathcal C}_G[\Gamma]$ induces a full
embedding $D^b({\mathcal C}_G[\Gamma'])\to D^b({\mathcal
C}_G[\Gamma])$ of bounded derived categories.  In particular,
given $M,N\in{\mathcal C}_G[\Gamma']$,
 $$\Ext^\bullet_{{\mathcal
 C}_G[\Gamma']}(M,N)\cong\Ext^\bullet_{{\mathcal
 C}_G[\Gamma]}(M,N)\cong\Ext^\bullet_G(M,N).$$

It follows that $\Ext^1_G(P_{\mu,<\lambda s},\text{\rm
rad}(\Delta(\lambda s)))=0$ and $\Ext^2_G(P_{\mu,<\lambda
s},\text{\rm rad}\Delta(\lambda s))=0$, so that
$\Ext^1_G(P_{\mu,<\lambda s},\Delta(\lambda
s))\cong\Ext^1_G(P_{\mu,<\lambda s},L(\lambda s))$.  Similarly, we
have the result (\ref{transferofextone}) already noted.
    So,  for any map
$P_{\mu,<\lambda s}\to\Delta(\lambda)$, there is a commutative
diagram
 \begin{equation}\label{commdiag}
 \begin{array}{ccc}
  \Ext^1_G(\Delta(\lambda),\Delta(\lambda s))
&
\begin{picture}(50,5) \put(0,3){\vector(1,0){50}}
\end{picture}
& \Ext^1_G(P_{\mu,<\lambda s},\Delta(\lambda s))
\\[2mm]
\begin{picture}(5,50) \put(0,50){\vector(0,-1){50}}
\end{picture}
  &&
\begin{picture}(5,50) \put(0,50){\vector(0,-1){50}}
\end{picture}
\\
\Ext^1_G(\Delta(\lambda),L(\lambda s)) &
\begin{picture}(50,5) \put(0,3){\vector(1,0){50}}
\end{picture}
& \Ext^1_G(P_{\mu,<\lambda s},L(\lambda s))
\end{array}
\end{equation}
in which the vertical maps are isomorphisms.

Now suppose that $P_{\mu,<\lambda
s}\overset\phi\to\Delta(\lambda)$ is not a split surjection.  Then
$\phi$ is not surjective, since $\Delta(\lambda)\in{\mathcal
C}_G[\Gamma']$ is projective. So $\phi$ factors as
$P_{\mu,<\lambda s}\to\text{\rm
rad}(\Delta(\lambda))\subseteq\Delta(\lambda)$.

Now assume that (\ref{Extcondition}) holds. Then the canonical
nonzero extension of $\Delta(\lambda)$ by $L(\lambda s)$ is zero
on $\text{\rm rad}(\Delta(\lambda))$, so the map
$\Ext^1_G(\Delta(\lambda),L(\lambda s))\to\Ext^1_G(P_{\mu,<\lambda
s},L(\lambda s))$ must be zero, as is the map
$\Ext^1_G(\Delta(\lambda),\Delta(\lambda
s))\to\Ext^1_G(P_{\mu,<\lambda s},\Delta(\lambda s))$.  Thus,
(\ref{Extcondition}) implies (\ref{theconditionstar}).

Conversely, suppose condition (\ref{theconditionstar})  holds. Let
$P_{<\lambda s}$ be a direct sum of sufficiently many
$P_{\mu,<\lambda s}$ with $\mu<\lambda$, $\mu\in\Gamma\cap
W_p\cdot 0$, so that there is a surjection $P_{<\lambda
s}\to\text{\rm rad}(\Delta(\lambda))$. Obviously, the component
maps $P_{\mu,<\lambda s}\to P_{<\lambda s}\to \text{\rm
rad}(\Delta(\lambda))$ do not give split surjections when composed
with the inclusion $\text{\rm
rad}(\Delta(\lambda))\subseteq\Delta(\lambda)$.  Thus, the
composites $\Ext^1_G(\Delta(\lambda),\Delta(\lambda
s))\to\Ext^1_G(P_{\mu,<\lambda s},\Delta(\lambda s))$ are zero by
(\ref{theconditionstar}). Using the above diagram, we obtain that
the map
$$\Ext^1_G(\Delta(\lambda),\Delta(\lambda s))\to\Ext^1_G(\text{\rm
rad}(\Delta(\lambda)),L(\mu))\to\Ext^1_G(P_{\mu,<\lambda
s},L(\lambda s))$$ is zero for each component $P_{\mu,<\lambda s}$
of $P_{<\lambda s}$. Thus, the map
$$\Ext^1_G(\Delta(\lambda),L(\lambda s))\to\Ext^1_G(\text{\rm
rad}(\Delta(\lambda)),L(\lambda s))\to\Ext^1_G(P_{<\lambda
s},L(\lambda s))$$
  is also zero.  However, the long exact sequence of $\Ext$ shows
  that that the map $\Ext^1_G({\text{\rm
  rad}}(\Delta(\lambda)),L(\lambda s))\to\Ext^1_G(P_{<\lambda
  s},L(\lambda s))$ is injective, since $\Hom_G(N,L(\lambda s))$
  for the kernel $N$ of $P_{<\lambda s}\to\text{\rm
  rad}(\Delta(\lambda))$.  This proves that the validity of
  (\ref{theconditionstar}) for all possible $\mu$, $\lambda,\lambda s$ as above,
  implies the validity of (\ref{Extcondition}).
\end{proof}

The following results ties the above equivalence in with the
even--odd varnishing equivalence (as a sufficient condition) for
the Lusztig conjecture, see, e.~g. \cite[\S7]{PS}.  See also
\S\ref{extoneisomorhismandhigherext} below.

\begin{prop} Let $\Gamma$ be a finite ideal in $X^+$.
Fix $\lambda,\lambda s\in \Gamma$ with $s\in \wSp$, such that
$\lambda<\lambda s$. Suppose that condition
(\ref{theconditionstar}) fails for some $\mu<\lambda$. If
$\mu<\lambda$ is chosen maximal for which this failure occurs,
then $\Ext^1_G(\Delta(\mu),L(\lambda))\not=0$ and
$\Ext^1_G(\Delta(\mu),L(\lambda s))\not =0$.

In particular, condition (\ref{theconditionstar}) holds if
$\Ext^1_G(\Delta(\mu),L(\gamma))=0$ always whenever
$\gamma\in\Gamma\cap W_p\cdot 0$ and $\mu,\gamma$ have lengths of
the same parity with respect to some ``length'' function
$l:\Gamma\cap W_p\cdot 0\to{\mathbb Z}$ for which $l(\lambda
s)\equiv l(\lambda)+1$ mod $2$.
 \end{prop}

\begin{proof}If (\ref{theconditionstar}) fails for some $\mu$, then the argument
for Thm. \ref{firstreduction} shows that the homomorphism
$\psi:\Ext^1_G(\Delta(\lambda),L(\lambda s))\to\Ext^1_G({\text{\rm
rad}}(\Delta(\lambda)),L(\lambda s))$ is not $0$.  Moreover, if
$E$ is a nonsplit extension $0\to L(\lambda s)\to E\to
\Delta(\lambda)\to 0$, and if $F$ is the largest submodule of $E$
without $L(\lambda s)$ as a composition factor, then $\psi$
factors through $\Ext^1_G(E',L(\lambda s))$, where $F'$ is the
image of $F$ in $\Delta(\lambda)$ and $E'=\Delta(\lambda))/F'$.
Also, (\ref{theconditionstar}) fails for $\mu$ if and only if
$L(\mu)$ is a composition factor of $E'=\text{\rm
rad}(\Delta(\lambda))/F'$. Taking $\mu$ maximal among such
weights, there are nonzero homomorphisms $\Delta(\mu)\to E'$ and
$\Delta(\mu)\to E^{\prime *}$, where $E^{\prime *}$ is the dual of
$E'$. (As well-known, the category ${\mathcal C}_G[\Gamma]$ has a
duality $M\mapsto M^*$ which fixes irreducible modules.) Pulling
back to $\Delta(\mu)$ the evident extensions of $E'$ (resp.,
$E^{\prime *}$) by $L(\lambda s)$ (resp., $L(\lambda)$) gives the
required nonsplit extensions of $\Delta(\mu)$ by $L(\lambda s)$
(resp., $L(\lambda)$).
\end{proof}

The papers \cite{CPS2}, \cite{CPS3} contain many other conditions
which are equivalent to the validity of the Lusztig character
formula. We will return to the even-odd vanishing result mentioned
in the proof of Thm. \ref{lusztigequiv} later in
\S\ref{extoneisomorhismandhigherext}.

\section{Type $A$}\label{typeA} We now consider the case $G=SL_n(k)$.
 {\it For technical reasons, we will assume throughout this
section that $p>3$. See the discussion below concerning
(\ref{contravariantequivalence}).}

For positive integers $n,r$, let (as before) $\Lambda^+(n,r)$ be
the set of all partitions of $r$ of length at most $n$. We regard
$\Lambda^+(n,r)$ as a poset, using the dominance ordering
$\trianglelefteq$. We will need to make use of a well-known method
of associating to any partition with at most $n$ nonzero parts a
dominant weight in $X^+$. For this, we will use a {\it local}
notational convention (i.~e., it will be used only in this
section!): We will denote partitions by symbols $\ulambda,\umu$,
etc.  Thus, if $\ulambda\in\Lambda^+(n,r)$, then
$\ulambda=(\ulambda_1,\ulambda_2,\cdots,\ulambda_n)$. Label the
fundamental dominant weights $\varpi_1,\cdots,\varpi_{n-1}$ for
$SL_n(k)$ as in \cite[p. 250]{Bo}. We have a mapping
$\Lambda^+(n,r)\to X^+$ which associates to
$\ulambda\in\Lambda^+(n,r)$ the dominant weight
$\lambda=\sum_{i=1}^{n-1}a_i\varpi_i$, where
$a_i=\ulambda_i-\ulambda_{i+1}$.

If a partition $\ulambda\in\Lambda^+(n,r)$ corresponds to a
dominant weight $\lambda\in X^+$ as above, then $\lambda$ is
regular (in the sense of \S\ref{charformulaholds}) if and only if
 \begin{equation}\label{regularcondition}
 \ulambda_i-\ulambda_{j}\not\equiv i-j\,\,{\text{\rm
mod}}\,\,p,\,\,\forall i<j\leq n\end{equation}
  (This alcove geometry notion of regularity is different than the Young
  diagram notion of $p$-regularity, which requires that no row be
  repeated $p$ times.)

 For a fixed pair $(n,r)$ of positive integers, let $S(n,r)$ be the
Schur algebra over $k$ of bidegree $(n,r)$.  Our notation will
largely be consistent with that in \cite[\S4]{CPS6}, where the
reader can find for more details on Schur algebras in the spirit
of this paper. (Also, in the notation of \S\ref{kunnethformula} in
Part I, $S(n,r)=S_q(n,r)$ for $q=1$.) The category $S(n,r)$--mod
identifies with the category of rational $GL_n(k)$-modules which
are polynomial of homogeneous degree $r$. The category
$S(n,r)$--mod is a highest weight category with poset
$(\Lambda^+(n,r),\trianglelefteq)$.

 Now fix a finite ideal $\Gamma\subset X^+$. Choose a sufficiently
large $r\equiv 0$ mod $n$ such that every weight in $\Gamma\cap
W_p\cdot 0$ corresponds (as described above) to a partition
$\ulambda\in \Lambda^+(n,r)$. Any $S(n,r)$-module $M$ is naturally
a $GL_n(k)$-module and hence an $SL_n(k)$-module, and ${\mathcal
C}_G[\Gamma\cap W_p\cdot 0]$ is a full subcategory of the image of
$S(n,r){\text{\rm--mod}}\to SL_n(k)$--mod under this
identification functor, with standard, costandard, and irreducible
modules going to standard, costandard, and irreducible modules,
respectively. In addition, if $M,N\in S(n,r)$--mod, then
\begin{equation}\label{identificationofext}
\Ext^\bullet_{S(n,r)}(M,N)\cong\Ext^\bullet_{GL_n(k)}(M,N)\cong
\Ext^\bullet_{SL_n(k)}(M,N),
\end{equation}
The first isomorphism is well-known (see \cite[(6)]{CPS6} and the
references there) and the second one follows immediately from an
elementary Hochschild-Serre spectral sequence argument, using the
normal subgroup $G=SL_n(k)$ of $GL_n(k)$.

Let mod--$k{\mathfrak S}_r$ denote the category of all finite
dimensional modules for the symmetric group ${\mathfrak S}_r$ of
degree $r$ over $k$. We will make use of Schur-Weyl duality between
$S(n,r)$--mod and mod--$k{\mathfrak S}_r$. In particular, let
$S(n,r)$--mod$(\Delta)$ be the full, exact subcategory of finite
dimensional $S(n,r)$-modules which have a $\Delta$-filtration,
i.~e., a filtration with sections of the form $\Delta(\ulambda)$,
$\ulambda\in\Lambda^+(n,r)$. Similarly, let mod--${\mathfrak
S}_r({\sS}_{n,r})$ be the full, exact subcategory of all finite
dimensional (right) $k{\mathfrak S}_r$-modules which have a
filtration with sections of the form $S_{\ulambda}$,
$\ulambda\in\Lambda^+(n,r)$, where $S_{\ulambda}$ is the Specht
module for $k{\mathfrak S}_r$ corresponding to the partition
$\ulambda$. Now assume that $p\geq n=h$ and recall the blanket
assumption that $p>3$. By \cite[(3.8.3.2)]{CPS5},\footnote{The
argument given in \cite{CPS5} is incomplete when $n=p$, unless
$n=r$. However, because of the assumption that $p>3$, the argument
is easily repaired using \cite[(3.3.4)]{HN} or the $i=1$ case of
\cite[Thm. 4.5]{PS}. In fact, these more recent results imply that
the assumption that $p\geq n$ may even be dropped in the equivalence
(\ref{contravariantequivalence}), when $p>3$. A related version of
(\ref{contravariantequivalence}) for $p>n$ is essentially a result
in \cite[p. 124]{E1}, restated in \cite[ftn. 15]{CPS5}. } there is a
contravariant equivalence
  \begin{equation}\label{contravariantequivalence}
S(n,r){\text{\rm --mod}}(\Delta)\overset{\text{\rm
contra}}{\underset\sim\longrightarrow}{\text{\rm
mod--}}k{\mathfrak S}_r(\sS_{n,r}).
\end{equation}
 If $\underline\mu\in\Lambda^+(n,r)$, let ${\mathfrak S}_{\underline\mu}$ be the
 associated Young (or parabolic) subgroup of ${\mathfrak S}_r$. Let
 $T_{\underline \mu}=\text{\rm Ind}_{{\mathfrak S}_{\underline\mu}}^{{\mathfrak S}_r}k$ be
 the corresponding permutation module obtained by inducing the trivial module $k$
 from ${\mathfrak S}_{\underline\mu}$ to ${\mathfrak S}_r$. If $V$
 is the standard $n$-dimensional module of column vectors for
 $GL_n(k)$, then ${\mathfrak S}_r$ acts on tensor space
 $V^{\otimes r}$ by place permutation.  As a right $k{\mathfrak
 S}_r$-module, $V^{\otimes r}$ decomposes into a direct sum of
 permutation modules $T_\ulambda$, $\ulambda\in\Lambda^+(n,r)$
 (each summand occurring with some positive multiplicity). Then
  $\Hom_{k{\mathfrak S}_r}(V^{\otimes r},V^{\otimes r})\cong
 S(n,r)$.
 The contravariant equivalence (\ref{contravariantequivalence}) is
given explicitly by restricting to the exact subcategory
$S(n,r){\text{\rm --mod}}(\Delta)$ the functor
\begin{equation}\label{diamondfunctor}
\begin{aligned}
S(n,r){\text{\rm --mod}}&\to{\text{\rm mod--}} k{\mathfrak S}_r,\\
& M\mapsto M^{\diamond}:=\Hom_{S(n,r)}(M,V^{\otimes r}).
\end{aligned}\end{equation} The equivalence (\ref{contravariantequivalence})
takes the standard module $\Delta(\ulambda)\in S(n,r){\text{\rm
--mod}}$ to the Specht module $S_\ulambda$. Also, the permutation
module $T_{\ulambda}$ corresponds to a projective
 module $P_\ulambda:=\Hom_{k{\mathfrak S}_r}(T_{\ulambda},V^{\otimes r})$
 of the type discussed in
 \S\ref{equivalentconditionssection}, {\it viz.,} it is a direct
 sum of the projective indecomposable module $P(\ulambda)$
  (with multiplicity one), together with various
 $P(\utau)$, for $\utau\triangleright\ulambda$, letting
 $\trianglerighteq$ denote the dominance order on
 $\Lambda^+(n,r)$. For more details, see \cite{CPS5},
 \cite{DPS1}, and \cite{DPS2}.

As a consequence of (\ref{contravariantequivalence}) and the
uniqueness of the short exact sequence
(\ref{shortexactsequence1}), we have the following result. (Recall
that $p>3$ throughout this section.)

\begin{lem}\label{spechtnonsplitextension}Let $\Gamma$ be a finite
ideal in $X^+$ as above. Choose $r\equiv 0$ mod$\,n$ so that every
$\lambda\in\Gamma\cap W_p\cdot 0$ corresponds to a partition
$\ulambda\in\Lambda^+(n,r)$. Suppose $\lambda<\lambda s$ both
belong to $\Gamma\cap W_p\cdot 0$ (where $s\in \wSp$).  Then, up
to scalar equivalence, there is a unique non-split short exact
sequence
\begin{equation}\label{spechtnonsplitextension1}
0\to S_\ulambda\to F\to S_{\underline{\lambda s}}\to 0
\end{equation}
in $k{\mathfrak S}_r$--mod.
\end{lem}

Each weight $\nu\in\Gamma\cap W_p\cdot 0$ is represented by a
partition $\underline{\nu}\in\Lambda^+(n,r)$.  The weight $0$ is
represented $\underline{0}=(r/n,\cdots, r/n)$.
  The partitions $\underline\nu$ with $\nu\in\Gamma\cap
W_p\cdot 0$ may be recursively defined as follows: Suppose this
notation is understood for all small partitions with at most $n$
parts.  Then $\nu\in\Gamma\cap W_p\cdot 0$ if and only if there is
a partition $\underline{\omega}\vartriangleleft\underline{\nu}$,
with ${\omega}\in \Gamma\cap W_p\cdot 0$, such that, for some
$0\leq i< j\leq n$ and $0<m$, we have:
\begin{itemize}
\item[(a)] $\unu_i=\uomega_i+m$ and $\unu_j=\uomega_j-m$;
 \item[(b)] $\unu_i-\unu_j+j-i\equiv m$ mod $p$.
\end{itemize}
The above conditions are just a translation of conditions for
$\omega$ to be related to $\nu$ by a certain kind of reflection.
Notice that $\unu_i-\unu_j+j-i< np<\uomega_i-\uomega_j+j-i$ for
some integer $n$. If the pair $i,j$ are unique with this property,
and only $m<p$ works in (a) above, then $\omega=\nu s$ (and
conversely, if $\uomega>\unu$).

Given $\umu\in\Lambda^+(n,r)$, $\mu\in\Jan$ if and only if
\begin{equation}\label{jantzenregionforpartitions}
\umu_1-\umu_n+n-1\leq p(p-n+2). \end{equation}
 As before, we
will consider if the Lusztig character formula holds for
$\Gamma\cap W_p\cdot 0$ for any ideal $\Gamma$ in $X^+$.

There is a natural and unique submodule
   $T_{\underline{\mu},<\underline{\gamma}}$, for any
   $\underline\gamma\in\Lambda^+(r)$, which has a filtration by Specht
   modules $S_{\underline\tau}$ with
   $\underline{\tau}\vartriangleleft\underline{\gamma}$, and with
   $T_{\underline \mu}/T_{\underline{\mu},<\underline{\gamma}}$
filtered by $S_{\underline\tau}$ with $\underline\tau$ not smaller
that $\underline\gamma$. The existence and
 uniqueness of $T_{\umu,<\ugamma}$ follow from
 (\ref{contravariantequivalence}) and the existence and uniqueness of the module
 $P_{\umu,<\ugamma}$ from \S\ref{equivalentconditionssection} (which follows from well-known
 quasi-hereditary algebra theory).\footnote{There is always a
 natural construction of $T_{\umu,<\ugamma}$, obtained by reduction mod $p$
  from a similar module
 defined over a principal ideal domain $Z$ of characteristic zero.
 See \cite[Thm. 5.2.1]{CPS5} for the existence of the relevant
 filtrations over $Z$, which uses no algebraic group theory.  The
 uniqueness of the analogue of $T_{\umu,<\ugamma}$ over $Z$ is
 easy. Indeed, over $Z$, $T_{\umu,<\ugamma}$ may be described as
 the intersection of the ($Z$-version) of $T_{\umu}$ with an
 evident canonical module defined over the quotient field of $Z$.}

 \begin{thm}\label{typeareduction}
Assume that $p\geq n$ and let $\Gamma$ be a finite ideal in $X^+$.
Then the Lusztig character formula holds for $\Gamma\cap W_p\cdot
0$ if and only if for each map $S_{\underline\lambda}\to
T_{\underline{\mu}, < \underline{\lambda s}}$, which is not a
split injection,  the induced map
\begin{equation}\label{stardiamond}
 \Ext^1_{k{\mathfrak S}_r}(S_{\underline{\lambda
 s}},S_{\underline{\lambda}})\to\Ext^1_{k{\mathfrak S}_r}(S_{\underline{\lambda
 s}},T_{\underline{\mu},<\underline{\lambda s}})
 \end{equation}
 on $\Ext^1$-groups
 is zero, whenever $\mu<\lambda<\lambda s$ in $\Gamma\cap W_p\cdot
 0$ (for some $s\in \wSp$). \end{thm}
\begin{proof} This follows from Thm. \ref{firstreduction} and
(\ref{identificationofext}) using the contravariant equivalence of
exact categories (\ref{contravariantequivalence}).
\end{proof}

Using Lemma \ref{spechtnonsplitextension}, we can restate the
above result in a slightly different, but suggestive way.

\begin{thm}Assume that $p\geq n$ and let $\Gamma$ be a finite ideal in $X^+$.
Then the Lusztig character formula holds for $\Gamma\cap W_p\cdot
0$ if and only if, whenever $\mu<\lambda<\lambda s$ belong to
$\Gamma\cap W_p\cdot 0$ (for some $s\in \wSp$), any morphism
$S_\ulambda\to T_{\umu,<{\underline{\lambda s}}}$ which is not a
split injection, fits into a commutative diagram
\begin{equation}
\begin{array}{ccccccccc}
0 &
\begin{picture}(17,0)
\put(0,5){\vector(1,0){23}}
\end{picture}
& S_\ulambda
\begin{picture}(0,0)
\put(-6,-6){\vector(0,-1){32}}
\end{picture}
&
\begin{picture}(20,0)
\put(-8,-32){\vector(-3,-2){10}}
\multiput(31,-6)(-13,-8.7){3}{\line(-3,-2){10}}
\put(-10,5){\vector(1,0){29}}
\end{picture}
& E &
\begin{picture}(21,0)
\put(-2,5){\vector(1,0){23}}
\end{picture}
& S_{\underline{\lambda s}}&
\begin{picture}(21,0)
\put(-2,5){\vector(1,0){23}}
\end{picture}
& 0
\\
&&&&&&&&
\\[7mm]
&& T_{\umu,<\underline{\lambda s}} &&&&&&
\end{array}
\end{equation}

\medskip\medskip
\medskip\noindent
of $k{\mathfrak S}_r$-modules, in which the horizonal row is the
unique (up to scalar equivalence) short exact sequence defined in
(\ref{spechtnonsplitextension1}).\end{thm}

Taking $\Gamma=\Jan$, the above theorems present a specific
necessary and sufficient condition in terms of symmetric group
modules for the validity of the Lusztig conjecture for $SL_n(k)$.
This result should be compared to \cite[Prop. 7.1]{PS}, where the
authors proved a sufficient condition for the validity of the
Lusztig conjecture for $SL_n(k)$. The condition involved only the
cohomology of symmetric groups, stating that the Lusztig
conjecture was true provided certain groups $\Ext^1_{k{\mathfrak
G}_s}(S_{\utau},D_{\usigma})$ were equal to 0. Here the (usually
large) integer $s$ is defined in terms of $n$ and the weights
$\Jan\cap W_p\cdot 0$. Also, $\utau,\usigma$ are various (regular)
partitions of $s$ indexed by the set $\Jan\cap W_p\cdot 0$, which
are defined using the Erdmann function $d$ in \cite[(6.1.2)]{PS}.
We also have the following result which is in the spirit of
\cite{PS}.

\begin{prop}\label{lusztigprop}Assume that $p\geq n$. Let $\Gamma\subset X^+$ be a finite
ideal as before.  Suppose, for fixed $\lambda<\lambda s$ as above,
that $(\ref{stardiamond})$ fails for some $\mu<\lambda$.  Let
$\mu<\lambda$ be maximal so that $(\ref{stardiamond})$ fails.
Assume that $\mu$, $\lambda$ and $\lambda s$ are all
$p$-restricted. Then $\Ext^1_{k{\mathfrak
S}_r}(S_{\umu'},D_{({\underline{\lambda s}})'})\not =0$ and
$\Ext^1_{k{\mathfrak S}_r}(S_{\umu'},D_{\ulambda'})\not=0$.
\end{prop}
\begin{proof}By \cite[Thm. (4.3)]{DEN} (or \cite[Cor.
5.3]{PS}),
\begin{equation}\label{den}
\dim\,\Ext^1_{S(r,r)}(\Delta(\ulambda),L(\umu))
\leq\dim\,\Ext^1_{k{\mathfrak S}_r}(S_{\ulambda'},D_{\umu'}),
\end{equation}
 whenever the partitions  $\ulambda'$ and $\umu'$ are $p$-regular.
For any partition $\ulambda$, the weight $\lambda\in X^+$ is
restricted if and only if the dual partition $\ulambda'$ is
$p$-regular. Now the result follows from (\ref{den}), together with
Thm. \ref{firstreduction} and
(\ref{identificationofext}).\end{proof}

\begin{rems}\label{erdmannremark} (a) Assume that $\Gamma=\Jan$ and that $p\geq 3p-2$ so
that all restricted weights lie in $\Jan$. By \cite{K} (see
 \S\ref{charformulaholds}), to see that the Lusztig character formula holds for
$\Jan\cap W_p\cdot 0$, it is enough to check it for restricted
regular weights. Then the condition that $\lambda, \lambda s$ be
restricted in the statement of Prop. \ref{lusztigprop} is not so
severe: Let $\tau\in X^+\cap W_p\cdot 0$, and suppose $\tau =
\lambda s$ for some $\lambda\in X^+$, $\lambda<\tau$. Further, if
$\tau$ is restricted, then $\lambda$ is also restricted. This
follows since $\lambda$ and $\tau$  are separated by only one
affine hyperplane $H_{\alpha,mp}$, and $H_{\alpha,mp}$ separates
$\tau$ from $0$. So $H_{\alpha,np}$  is not among those
hyperplanes defining the restricted region. If the Lusztig
character formula holds for each weight in $X^+\cap W_p.0$ smaller
than
 $\tau$, then Thm. \ref{lusztigequiv}(b)  shows that it holds for
$\tau$ if and only the inequality (\ref{Extcondition}) holds.
Thus, in Prop. \ref{lusztigprop}, we can assume to start that both
$\lambda$ and $\lambda s$ restricted. (However, it still may be
true that $\mu$ is not restricted. The Erdmann function $d$
\cite[(6.1.2)]{PS} in offers some way around this problem, but it
would seem best to first check if the conclusion of Prop.
\ref{lusztigprop} holds when all three weights involved are
restricted.)

(b) The results of this section avoid the use of the
 function $d$ altogether. Thus, they
 give a way of using a smaller symmetric group than in Erdmann's paper \cite{E2} to
 determine, say, the decomposition  numbers for the ``principal"
block
 (containing the determinant module) of the Schur algebra $S(n,r)$
 in terms of ${\mathfrak S}_r$, without increasing $r$.
 However, it has been long
 known that the decomposition numbers of the endomorphism ring of a
 permutation module can be equated
 with multiplicities of ordinary characters inside those of $p$-adic
 indecomposable components of the
 permutation module. See \cite{S1}. So, Schur algebra decomposition
 numbers can be computed in terms of symmetric group permutation module theory
in
 this way, though this method would require knowledge
the indecomposable components of permutation modules.  Erdmann's
reduction notably uses (albeit after a big increase in $r$) just
the projective
 indecomposable components, whose multiplicities are actual symmetric group
 decomposition numbers. Perhaps the notable aspect of our reduction, aside
 from the fact that no increase in $r$ is required, is that no indecomposable
 components of permutation modules are required at all, only certain
 submodules of permutation modules whose definition is transparent in terms of
 ordinary character theory.
 \end{rems}

\medskip\medskip
\begin{center} {\large {\bf PART III: Levi subgroups and irreducible modules}}\end{center}
\addcontentsline{toc}{section}{\bf Part III: Levi subgroups and
irreducible modules }
\medskip

Again, $G$ is a simple, simply connected algebraic group over $k$.
When $p\geq h$, the discussion will often make use of the notion,
defined in \S\ref{charformulaholds}, of what it means for the
Lusztig character formula to hold for a subset $\Sigma \subseteq
W_p\cdot 0$.
\section{A theorem of Hemmer and a generalization}\label{hemmersection}
It is natural to ask for a version of (\ref{transfer1}) relating
the groups $\Ext^\bullet_{G}(L(\lambda),L(\mu))$ to the groups
$\Ext^\bullet_H(L_H(\lambda),L_H(\mu))$. This section considers
this question. Thm. \ref{Hemmer} below establishes an inequality
for $\Ext^1$-groups which is in the spirit of (\ref{transfer1})
and (\ref{transfer2}). A version of this result for quantum
enveloping algebras at a root of unity is easily obtained along
the same lines.

In \cite[Thm. 2.3]{H}, Hemmer proves that if
$\lambda,\mu\in\Lambda^+(n,r)$ satisfy $\lambda_1=\mu_1=m$ then
there is an injection
\begin{equation}\label{Hemmerinjection}
\Ext^1_{S(n,r)}(L(\lambda),L(\mu))\hookrightarrow\Ext^1_{S(n-1,r-m)}(L(\bar\lambda)),L
(\bar\mu)),
\end{equation}
  where $\bar\lambda$ is the partition of $r-m$ obtained
by removing $\lambda_1$ from $\lambda$, and $\bar\mu$ is defined
similarly. This result is then interpreted in terms of comparing
the $\Ext$-quiver of $S(n,r)$ with that of $S(n-1,r-m)$.

 We can
generalize Hemmer's theorem to other types:

 Let  $H$ be a Levi subgroup of $G$ and
$\Omega=\omega+{\mathbb Z}\Phi_H$ is a coset of ${\mathbb
Z}\Phi_H$ in ${\mathbb Z}\Phi$. Put $\Omega^+=\Omega\cap X^+$.
\begin{thm}\label{Hemmer}Let $\lambda,\mu\in \Omega^+$.
 Then there is a
natural injection
\begin{equation}\label{injection}
\Ext^1_G(L(\lambda),L(\mu))
\hookrightarrow\Ext^1_H(L_H(\lambda),L_H(\mu)).\end{equation}
\end{thm}

\begin{proof}We will argue that the required injection is
induced by the truncation functor $\pi_\Omega$
(\ref{truncationfunctor}).  We can assume that $\lambda\leq\mu$.
There is a natural short exact sequence $0\to
Q(\lambda)\to\Delta(\lambda)\to L(\lambda)\to 0$. The long exact
sequence of $\Ext^\bullet_G$ with respect to the functor
$\Hom_G(-,L(\mu))$ gives an exact sequence
$$\Hom_G(Q(\lambda),L(\mu))\to\Ext^1_G(L(\lambda),L(\mu))\to
\Ext^1_G(\Delta(\lambda),L(\mu)).$$
  All the composition factors $L(\tau)$ of
$Q(\lambda)$ satisfy
  $\tau<\lambda<\mu$; so
$\Hom_G(Q(\lambda),L(\mu))=0$.
  A similar conclusion holds for the analogous exact sequence
  $0\to Q_H(\lambda)\to \Delta_H(\lambda)\to L_H(\lambda)\to 0$.
  Therefore, the functor $\pi_\Omega$ provides a
commutative diagram
$$
\begin{array}{ccccccccc}
0 \hspace*{6pt} & \to & \Ext^1_G(L(\lambda),L(\mu))
   & \to & \Ext^1_G(\Delta(\lambda),L(\mu))
  \\[2mm]
&&
\begin{picture}(10,50) \put(8,50){\vector(0,-1){50}}
\put(-4,25){$\scriptstyle{\pi_\Omega}$}
\end{picture}
&&
\begin{picture}(10,50) \put(8,50){\vector(0,-1){50}}
\put(-4,25){$\scriptstyle{\pi_\Omega}$}
\end{picture}
&& \\[2mm]
0 \hspace*{6pt} & \to & \Ext^1_H(L_H(\lambda),L_H(\mu))& \to &
\Ext^1_H(\Delta_H(\lambda),L_H(\mu))
\end{array}
 $$
  with exact rows.  By (\ref{Extisos}) for $M=L(\mu)$ (or \cite[Cor. 10]{CPS6}), the right hand
vertical map is an isomorphism.  Hence, the left hand vertical map
is an injection, as required. \end{proof}

\section{An example}\label{example} In general, the inclusion map
$\Ext^1_G(L(\lambda),L(\mu))\hookrightarrow\Ext^1_H(L_H(\lambda),L_H(\mu))$
given in Thm. \ref{Hemmer} need not be an
isomorphism.\footnote{The examples in this section provide a
(negative) answer to a question raised by Hemmer in his lecture at
the AMS Conference ``Representations of algebraic groups, quantum
groups, and Lie algebras," Snowbird, Utah, July 2004. The form of
the examples suggest a positive answer to Hemmer's question is
almost equivalent to the validity of the Lusztig conjecture.}

Let $G=SL_3(k)$, with char$\,k=p=3$.  We list the dominant weights
as in \cite[p. 250]{Bo}. The Lusztig conjecture is known to be
true for $G$. In particular, by the result of Kato discussed in
\S\ref{charformulaholds}, the Lusztig character formula holds for
$4\varpi_1+\varpi_2= (\varpi_1+\varpi_2) + 3\varpi_1$, since
$\varpi_1+\varpi_2$ is a restricted weight in $W_p\cdot 0$ and
$\varpi_1\in\overline{C^+}$. Similarly, it holds for $3\varpi_1$.
Thus, the Lusztig character formula holds  for each weight in the
sequence $0<\varpi_1\varpi_2< 3\varpi_1<4\varpi_1+\varpi_2$ of
adjacent dominant weights. On the other hand, Lemma
\ref{amazinglemma} (or an unpleasant computation) shows that the
Lusztig character formula fails for
$3\varpi_1+3\varpi_2=(4\varpi_1+\varpi_2)s_{\alpha_2}$.  By Thm.
\ref{lusztigequiv}(c),
$\Ext^1_G(L(4\varpi_1+\varpi_2),L(3\varpi_1+3\varpi_2))=0$.
However, let $H$ be the Levi subgroup with $\Phi_H=\{\alpha_2\}$
(i.~e., defined by $\Pi_1$ in the notation of
\S\ref{kunnethformula}), then, directly or using Thm.
\ref{lusztigequiv}(c) again, we have that
$\Ext^1_H(L_H(4\varpi_1+\varpi_2),L_H(3\varpi_1+3\varpi_2))
\not=0$.

In terms of the group $GL_3(k)$, let
$\lambda=(6,2,1),\mu=(6,3,0)\in\Lambda^+(3,9)$. In the notation of
\S\ref{kunnethformula}, Example 1,
$$\Ext^1_{GL_3(k)}(L(\lambda),L(\mu)=0,$$
while
$$\Ext^1_H(L(\overline\lambda),L_H(\overline\mu))\not=0.$$
Therefore, the map (\ref{Hemmerinjection}) is not an isomorphism
in this case.

Another, similar but interesting example, is provided by $GL_5(k)$
when $p=5$, taking $\lambda=(10,5^3), (10,5^2,4,1)\in
\Lambda^+(5,25)$.
 (From computer results
 discussed in \cite{S2}, and the translation principle, (\ref{lusztigcharacterformula})
 holds in the category ${\mathcal C}_G$  for   any
 regular restricted weight $\nu$.)
In fact, this last example let us to the ``smaller" one above.

\medskip
Despite these examples,  \S\ref{extoneisomorhismandhigherext}
below shows that (\ref{injection}) {\it is} an isomorphism for
regular weights lying in a finite ideal $\Gamma$ with the property
that the Lusztig character formula holds for $\Gamma\cap W_p\cdot
0$.

\section{An $\Ext^1$-isomorphism result and higher
$\Ext^\bullet$-groups}\label{extoneisomorhismandhigherext}
  In this section, we take up, for arbitrary cohomological degree
  $n$,  the study of the map
  $\Ext^n_G(L(\lambda),L(\mu))
  \to\Ext^n_H(L_H(\lambda),L_H(\mu))$
  induced by the truncation functor $\pi_\Omega$.
We will show that if $\lambda,\mu$ are regular weights, lying in a
finite ideal $\Gamma$ such that the Lusztig character formula
holds for $\Gamma_{\text{\rm reg}}$, then $\pi_\Omega$ always
induces an {\it surjection} provided $\lambda,\mu$ belong to the
same ${\mathbb Z}\Phi_H$-coset. In particular, (\ref{injection})
is an isomorphism in this case. To obtain such a result, it will
be necessary to work in the context of the homological dual of a
quasi-hereditary algebra, and so we begin by recalling some
general results.

In \cite{CPS3}, the homological dual $A^!$ of a quasi-hereditary
algebra $A$ was defined to be the Yoneda $\Ext$-algebra
$$A^!=\Ext^\bullet_A(L_0,L_0)$$
where $L_0$ is the direct sum of the distinct irreducible
$A$-modules. If $\sC=A$--mod, write $\sC^!$ for $A^!$--mod.

 The rest of this section will make use of properties of
a highest weight category $\sC\cong A$--mod having a
Kazhdan-Lusztig theory with respect to a length function
$l:\Lambda\to{\mathbb Z}^+$ on its weight poset $\Lambda $. This
means that, for any $\lambda,\mu\in\Lambda$,
 \begin{equation}\label{KLtheory}
 \begin{cases}
 \Ext^n_{\sC}(L(\lambda),\nabla(\mu))\not=0 \\
 \,\,\,\,\,\,\,\,\,\, \text{\rm or}\\
\Ext^n_\sC(\Delta(\mu),L(\lambda))\not=0 \end{cases} \implies
n\equiv
 l(\lambda)-l(\mu)\,{\text{\rm mod}}\,\,2
\end{equation}
 See \cite[\S3]{CPS2} or \cite[\S1.3]{CPS3} for further discussion.
We will always assume below that $\sC$ has finite poset $\Lambda$,
so that $\sC\cong A$--mod for some finite dimensional algebra $A$.
The following result, proved in \cite{CPS3}, summarizes some basic
properties of the homological dual.

\begin{thm}\label{recoll} Assume that $\sC\cong A{\text{\rm
--mod}}$ has a Kazhdan-Lusztig theory with respect to a length
function $l:\Lambda\to{\mathbb Z}^+$.  Let $\Xi$ be a coideal in
$\Lambda$. Then:

(a) $A^!$ is a quasi-hereditary algebra; more precisely,
$\sC^!=A^!$--mod is a highest weight category with weight poset
$\Lambda^{\text{op}}$ (the poset {\it opposite} to $\Lambda$).

(b) $\sC(\Xi)$ has a Kazhdan-Lusztig theory with respect to
$l|_\Xi$.

(c) ${\mathcal C}(\Xi)^!\cong{\mathcal C}^![\Xi^{\text{\rm op}}]$
as highest weight categories.
\end{thm}

\begin{proof} (a) is proved in \cite[Thm. 2.1]{CPS3}.  (b) is
well-known, but follows directly in the spirit of this paper as
follows: Let $j^*:\sC\to\sC(\Xi)$ be the natural quotient functor.
By \cite[Lemma 6]{CPS6}, given $\lambda,\mu\in\Xi$,
$\Ext^\bullet_{\sC(\Xi)}(L_\Xi(\lambda),\nabla_\Xi(\mu))\cong
\Ext^n_{\sC}(L(\lambda),\nabla(\mu))$ and
$\Ext^\bullet_{\sC(\Xi)}(\Delta_\Xi(\mu),L_\Xi(\lambda)) \cong
\Ext^\bullet_\sC(\Delta(\mu),L(\lambda))$, where we have written
$L_\Xi(\lambda)$, etc. for the irreducible, etc. object in
$\sC(\Xi)$ indexed by $\lambda$. Now (b) follows from the
definition (\ref{KLtheory}). For (c), see \cite[(2.3)]{CPS3}.
\end{proof}

In (c) above, ${\mathcal C} (\Xi )^!\cong (A_{\Xi})^!$--mod, where
$A_{\Xi}$ is a quasi-hereditary algebra such that
$A_{\Xi}{\text{\rm --mod}} \cong{\mathcal C}(\Xi)$. The proof
given in \cite[p. 308]{CPS3} is more precise, showing that the
natural algebra homomorphism
 \begin{equation}\label{defnphi}\Phi:A^!=\Ext^\bullet_{{\mathcal
C}}(L_0,L_0)\to\Ext^\bullet_{{\mathcal C}(\Xi)}(j^*L_0,j^*L_0)
=(A_{\Xi})^!,\end{equation}
  induced by the quotient functor $j^*:{\mathcal
C}\to{\mathcal
  C}(\Xi)$,  defines  an
  equivalence
$$\Phi^*:(A_{\Xi})^!{\text{\rm--mod}}\to A^!/J{\text{\rm --mod}}$$
of highest weight categories, for some (idempotent) ideal $J$
(denoted $A^!s_\Gamma A^!$ in \cite{CPS3}) which is part of a
defining sequence for $A^!$. Thus, $A^!/J\cong (A_{\Xi})^!$.   In
particular, we have:

  \begin{cor} \label{surjectivephicor} Assume that $A{\text{\rm --mod}}$ has a
Kazhdan-Lusztig theory. Then the map $\Phi$ defined in
(\ref{defnphi}) is surjective.
  \end{cor}

We now return to the setting of a simple, simply connected
algebraic group $G$ over $k$, assuming that $p\geq h$. Fix a
finite ideal $\Gamma\subseteq X^+$, and let $\sC_{G,0}[\Gamma]$
denote the highest weight category consisting of all finite
dimensional $G$-modules with composition factors $L(\xi)$ with
$\xi\in\Gamma\cap W_p\cdot 0$.\footnote{If we replace the $\leq$
partial ordering on $X^+$ by the Jantzen $\uparrow$ ordering, then
$\Gamma\cap W_p\cdot 0$ becomes an ideal in $(X^+,\uparrow)$.}
Define a length function $l:\Gamma\cap W_p\cdot 0\to{\mathbb N}$
by $l(w\cdot 0)=l(w)$. By the argument in \cite[Thm. 5.3]{CPS2}
(see the discussion of Thm. \ref{lusztigequiv}),
$\sC_{G,0}[\Gamma]$ has a Kazhdan-Lusztig theory if and only if
 the Lusztig character formula holds for
 $\Gamma\cap W_p\cdot 0$.

Let $H$ be a Levi subgroup of $G$ as in \S\ref{sectionone}. In
Thm. \ref{maintheorem}, $\Gamma^+_F$ can be replaced by a ``block
analog". That is, let $\mathcal B$ be any union of $W_p$-orbits in
the dot action on $X$. Put $\Gamma_{F,{\mathcal
B}}^+=\Gamma^+_F\cap {\mathcal B}$. Then clearly
$\Omega^+_{F,{\mathcal B}}:=\Gamma^+_{F,{\mathcal B}}\cap \Omega$
is a union of dominant weights in orbits of the affine Weyl group
$W_{H,p}$ of $H$. (Use the same weight $\rho$ in defining the dot
action of $W_{H,p}$ as is used for $W_p$.) Thm. \ref{maintheorem}
remains true with $\Gamma_F^+$ and $\Omega_F^+$ replaced by
$\Gamma^+_{F,{\mathcal B}}$ and $\Omega^+_{F,{\mathcal B}}$,
respectively.

We will take ${\mathcal B}=X^+_{\text{\rm reg}}$, the set of
regular dominant weights, writing $\Gamma^+_{F,{\text{reg}}}$ and
$\Omega^+_{F,{\text{reg}}}$ for $\Gamma^+_{F,{\mathcal B}}$ and
$\Omega^+_{F,{\mathcal B}}$, respectively,

Fix a finite ideal $\Sigma$ in $X^+$ such that the Lusztig
character formula holds for $\Sigma_{\text{\rm
reg}}\not=\emptyset$. Let $\Omega=\omega+{\mathbb Z}\Phi_H$ be a
coset of ${\mathbb Z}\Phi_H$ in $X$, with $\omega\in\Sigma$. Put
$F=\Sigma\cap\Omega$ and, as in Thm. \ref{maintheorem}, form
$\Gamma^+_{F}$ and $\Omega^+_F$.  Thus, $\Omega^+_F$ is a coideal
in $\Gamma^+_F$ and an ideal in $(X^{(H)+},\leq)$. Set
$\Lambda_0=\Gamma^+_{F,{\text{\rm reg}}}$ and
$\Omega_0=\Omega^+_{F,{\text{\rm reg}}}$, and form the highest
weight categories ${\mathcal C}_{G,0}={\mathcal C}_G[\Lambda_0]$
and ${\mathcal C}_{H,0}={\mathcal C}_H[\Omega_0]$. Let $L_0$  be
the direct sum of the distinct irreducible $G$-modules
$L(\lambda)$ for $\lambda\in\Lambda_0$.  Similarly, let $L_{0,H}$
be the direct sum of the distinct irreducible $H$-modules
$L_H(\lambda)$ for $\lambda\in\Omega_0$.

 \begin{thm}\label{kltheorem}Assume the notation of the previous
 paragraph.
    Then:

  (a)  The category $\sC_{H,0}$ has a
Kazhdan-Lusztig theory.
  In particular, the algebra
$A^!_H:=\Ext^\bullet_{H}(L_{0,H},L_{0,H})$
  is quasi-hereditary.

  (b) The natural map
  $$A^!=\Ext^\bullet_G(L_0,L_0)\to\Ext^\bullet_H(L_{0,H},L_{0,H})=A_H^!$$
  induced by $\pi_\Omega$ is surjective, and its kernel is
a
  defining ideal $J$ of the quasi-hereditary algebra
  $A^!$.

  (c) We have $A^!/J\text{\rm
--mod}\cong\sC^!_{G,0}[\Omega_0^{\text{\rm
  op}}]$.
  \end{thm}

  \begin{proof}  (a) follows from Thm. \ref{recoll}, using
 Thm. \ref{maintheorem} which shows that
${\mathcal C}_{H,0}\cong{\mathcal C}_{G,0}(\Omega_0)$. Then (b)
and (c) follow from Cor. \ref{surjectivephicor}.
  \end{proof}

  \begin{rems}\label{ffinalremark} (a)A somewhat less precise way to state part
(c) is, in
  view of part (b),
$$\sC_{G,0}^![\Omega^{\text{\rm
op}}_0]\cong\sC_H[\Omega_0]^!.$$

(b) As remarked earlier, for semisimple algebraic groups having a
fixed root system, the Lusztig conjecture is known to hold
provided that the characteristic $p$ is large enough, though lower
bound is known on the size of $p$ in general. Thus, in those
cases, we can in Thm. \ref{kltheorem} take $\Sigma=\Jan$, the
Jantzen region.
\end{rems}

\begin{cor}\label{lusztigcor}Assume the set-up in the paragraph
immediately before Thm. \ref{kltheorem}. Let $\lambda,\mu\in
\Sigma_{\text{\rm reg}}\cap\omega+{\mathbb Z}\Phi_H$. Then

(a)
$\Ext^1_G(L(\lambda),L(\mu))\cong\Ext^1_H(L_H(\lambda),L_H(\mu)).$

(b) The map
$\Ext^n_G(L(\lambda),L(\mu))\to\Ext^n_H(L_H(\lambda),L_H(\mu))$
which is induced by $\pi_\Omega$ is surjective for all $n\geq 0$.
In particular,
$$\dim\,\Ext^n_G(L(\lambda),L(\mu))\geq\dim\,\Ext^n_H(L_H(\lambda),L_H(\mu)).$$
\end{cor}
\begin{proof}Part (b) follows directly from Thm.
\ref{kltheorem}. Part (a) is then a consequence of (b) for $n=1$
and Thm. \ref{Hemmer}.
\end{proof}

\begin{rems}(a) An alternative proof of the inequality in Cor.
\ref{lusztigcor}(b) can be based on the $\Ext$-formula proved in
\cite[Thm. 3.5]{CPS2}, which states that
$$\dim\Ext_{\sC}^n(L(\lambda),L(\mu))=\sum_{m=0}^n\sum_\nu
\dim\Ext_{\sC}^m(L(\lambda),\nabla(\nu))
\cdot\dim\Ext_{\sC}^{n-m}(\Delta(\nu),L(\mu))$$
 in any highest weight category $\sC$ with a Kazhdan-Lusztig theory.
This formula shows just what part of
$\dim\Ext^n_G(L(\lambda),L(\mu))$ contributes of
$\Ext^n_H(L_H(\lambda),L_H(\mu))$---namely, the weight $\nu$ must
 lie in $\Omega$.

(b)
 Now we consider Hemmer's
injective map (\ref{Hemmerinjection}) in the case when the Lusztig
conjecture is true for $SL_n$. Assume that
$\lambda,\mu\in\Lambda^+(n,r)$ satisfy $\lambda_1=\mu_1=m$, and
let $\bar\lambda,\bar\mu$ be as before. Assume that
$\lambda_i-\lambda_j\not\equiv i-j$ and $\mu_i-\mu_j\not\equiv
i-j$ mod$\,p$ for all $i<j\leq n$ (see (\ref{regularcondition}))
and that $ \lambda_1-\lambda_n+n-1\leq p(p-n+2)$ and
$\mu_1-\mu_n+n-1\leq p(p-n+2)$ (see
(\ref{jantzenregionforpartitions})). Then {\it the map
(\ref{Hemmerinjection}) is an
  isomorphism.} Moreover, for higher degree $\Ext$-groups,
there is
  a corresponding {\it surjection}, a somewhat surprising
  turnabout.

 (c) Observe that an alternative proof of Cor. \ref{lusztigcor}(a)
 can be based on the commutative diagram in the proof of Thm. \ref{Hemmer}, once it is
 observed that \cite[Thm. 4.3]{CPS2} implies that the
 horizonal map is surjective.
\end{rems}
 \begin{cor}\label{lastcor}Let $\Sigma_{\text{\rm reg}}$ be as in
the set-up in the paragraph immediately before Thm.
\ref{kltheorem}, and put ${\bold X}:=\Sigma_{\text{\rm reg}}$. (In
particular, the Lusztig character formula holds for irreducible
modules $L(\lambda)$, $\lambda\in{\bold X}$.) Let $\sC_H[{\bold
 X}]$ denote the category of all finite dimensional rational
 $H$-modules which have composition factors $L_H(\xi)$, $\xi\in\bold X$.
 Then $\sC_H[{\bold X}]$ is a highest weight category with respect
 to the poset $({\bold X},\leq_H)$, and, further, $\sC_H[{\bold X}]$
 has a Kazhdan-Lusztig theory in the sense of (\ref{KLtheory}).\end{cor}
 \begin{proof}First note, that if $\Delta_H(\nu)$ is a standard module for
 $\nu\in{\bold X}$, and $L_H(\varpi)$ is a composition factor, then $\varpi\in X^+$
since $X^+$ is an ideal in $(X^{(H)},+,\leq_H)$.  Also,
$L(\omega)$ is a composition
 factor of $\Delta(\nu)$---this is a result of Donkin
 \cite{D1}---see also \cite[Cor. 12]{CPS6}. In particular, $\nu$ and $\varpi$ belong
 to the same block of $G$, and so $\nu$ and $\varpi$ are conjugate under the dot action
 of $W_p$.  Hence, $\varpi\in\bold X$.
  Now consider the projective cover
 $P_{\Gamma,H}(\varpi)$ of
$L_H(\varpi)$ in the category $\sC_H[\Gamma]$, where $\Gamma$ is
the set of dominant weights (for $G$) regarded as an ideal in
$(X^{(H)+},\leq_H)$.  If $\Delta_H(\xi)$ is a $\Delta_H$-section,
then $L_H(\varpi)$ is a composition factor of $\nabla_H(\xi)$ by
Brauer-Humphreys reciprocity \cite[Thm. 3.11]{CPS1}.  Hence,
$L_H(\varpi)$ is a composition factor of $\Delta_H(\xi)$.  As
before, this forces $\xi\in\bold X$. Thus, $\sC_H[{\bold X}]$ is a
highest weight category.

   The final assertion follows from Thm.
 \ref{kltheorem}(a), since $\Ext_H^\bullet$ vanishes between
 modules with respective composition factors in different cosets
 of the root lattice of $H$.
 \end{proof}

Finally, we have the following general (and elementary) result
which shows that the $\Ext^1$-isomorphism property behaves well
with respect to Frobenius twisting.

\begin{thm}Assume that if $G$ or $H$ has a component of type
$C_m$, then $p>2$. Suppose that $\lambda,\mu$ lie in the same
${\mathbb Z}\Phi$-coset $\Omega=\omega+{\mathbb Z}\Phi_H$ and that
$\pi_\Omega$ induces an isomorphism
$\Ext^1_G(L(\lambda),L(\mu))\overset\sim\to\Ext^1_H(L_H(\lambda),L_H(\mu))$.
For any positive integer $r$, put $\Omega^{(r)}=p^r\omega
+{\mathbb Z}\Phi_H$. Then $\pi_{\Omega^{(r)}}$ induces an
isomorphism
$\Ext^1_G(L(p^r\lambda),L(p^r\mu))\overset\sim\to\Ext^1_H(L_H(p^r\lambda),L_H(p^r\mu))$.
\end{thm}
\begin{proof}It suffices to consider the $r=1$ case. The hypothesis
on $p$ guarantees that, if $G_1$ denotes the first Frobenius
kernel of $G$, then $H^1(G_1,k)=0$ \cite[Prop. 12.9]{Jan}. Hence,
by a Hochschild-Serre spectral sequence argument for $G_1\lhd G$,
$\Ext^1_G(L(p\lambda),L(p\mu))\cong \Ext^1_G(L(\lambda),L(\mu))$.
Similar remarks apply to $\Ext^1_H(L_H(p\lambda),L_H(p\mu))$.
Therefore,
$$\dim\,\Ext^1_G(L(p\lambda),L(p\mu))=\dim\,\Ext^1_H(L_H(p\lambda),L_H(p\mu))$$
and the result follows from Thm. \ref{Hemmer}, since
$p\lambda,p\mu\in\Omega^{(1)}$.
\end{proof}

Thus, in general, the injection (\ref{injection}) is an
isomorphism for some regular weights well outside $\Jan$, where
the Lusztig character formula fails.\footnote{Assume $p\geq h$.
Lemma \ref{amazinglemma} provides a way to construct $\lambda\in
W_p\cdot 0\cap X^+$ such that the Lusztig character formula holds
for $\lambda$ but fails for $p\lambda$.}

\end{document}